\documentclass{amsproc}
\usepackage{eurosym}
\usepackage{amssymb}
\usepackage{amsfonts}
\usepackage{cmtt}

\DeclareMathOperator{\Conn}{Conn}

\DeclareMathOperator{\Spec}{Spec}
\theoremstyle{plain}
\newtheorem{theorem}{Theorem}[section]

\newtheorem{claim}[theorem]{Claim}

\newtheorem{definition}[theorem]{Definition}
\newtheorem{example}[theorem]{Example}
\newtheorem{lemma}[theorem]{Lemma}
\newtheorem{proposition}[theorem]{Proposition}
\newtheorem{remark}[theorem]{Remark}

\newtheorem{innercustomlemma}{Lemma}
\newenvironment{customlemma}[1]
{\renewcommand\theinnercustomlemma{#1}\innercustomlemma}
{\endinnercustomlemma}

\theoremstyle{remark}
\numberwithin{equation}{section}
\newcommand{\field}[1]{\mathbb{#1}}
\newcommand{\C}{\field{C}}

\newcommand{\N}{\field{N}}

\newcommand{\R}{\field{R}}
\newcommand{\Z}{\field{Z}}

\newcommand{\superscript}[1]{\ensuremath{^{\textrm{#1}}}}

\def\wu{\superscript{*}}
\def\wg{\superscript{\dag}}

\begin{document}

\title{Gorenstein matroids}

\author[Micha{\l} Laso{\'n}]{Micha{\l} Laso{\'n}\wu\footnote{\wu michalason@gmail.com; Institute of Mathematics of the Polish Academy of Sciences, ul.\'{S}niadeckich 8, 00-656 Warszawa, Poland}}

\author[Mateusz Micha{\l}ek]{Mateusz Micha{\l}ek\wg\footnote{\wg mateusz.michalek@uni-konstanz.de; University of Konstanz, Fach D 197, 78457 Konstanz, Germany}}

\thanks{Research supported by Polish National Science Centre grant no. 2015/19/D/ST1/01180.}


\begin{abstract}
Matroid base polytope and independence polytope are Cohen--Macaulay. There are plenty of examples when they are also Gorenstein. We provide a full classification of matroids whose independence polytope or base polytope is Gorenstein. This answers to a question raised by Herzog and Hibi.
\end{abstract}

\maketitle

\section{Introduction}

Matroids and lattice polytopes are combinatorial objects that are fundamental for combinatorial algebraic geometry. They belong to a part of mathematics where the interaction of algebra and geometry with combinatorics is particularly strong and significant. We explore this connection.

\subsection{Algebraic motivation}

Toric varieties are a class of algebraic varieties that on the one hand capture many varieties seen in applications and on the other hand are more amenable to combinatorial techniques than general algebraic varieties. Indeed, the geometry of a toric variety is fully determined by the combinatorics of its associated lattice polytope. When an algebraic variety is constructed using only combinatorial data, one expects to have a combinatorial description of its algebraic properties. An attempt to achieve this description often leads to surprisingly deep combinatorial questions.

The toric variety of a matroid is a particularly interesting example. For a representable matroid it has a nice geometric description -- it is the torus orbit closure in a Grassmannian, moreover every orbit closure arises in this way \cite{GeGoMaSe87}. The toric variety of a matroid is recognized mostly due to a famous conjecture of White \cite{Wh80}. The conjecture provides a description of generators of the toric ideal of a matroid, that is the ideal defining the matroid variety. In particular, it states that this ideal is generated in degree two. The conjecture was confirmed for several special classes of matroids. For general matroids it was proved `up to saturation' \cite{LaMi14}, and later upgraded for `high degrees w.r.t. the rank' \cite{La16}. In full generality White's conjecture remains open since its formulation in $1980$. 

In this paper we study algebraic properties of the toric variety of a matroid. White \cite{Wh77} already proved that the toric variety of a matroid is normal. Hence, by a celebrated result of Hochster, it satisfies Cohen--Macaulay property. We investigate a natural stronger property, namely the Gorenstein property. The notion of the Gorenstein property goes back to Grothendieck. It reflects many symmetries of cohomological properties. It also implies that singularities of the variety are not `too bad'. In particular, there is the following chain of inclusions.
\begin{center} 
smooth varieties $\subset$ Gorenstein varieties $\subset$ Cohen--Macaulay varieties 
\end{center}

Though, not every matroid variety is Gorenstein. Thus, there is a need for a classification. This line of research was pioneered by Herzog and Hibi \cite{HeHi02} who classified `generic' discrete polymatroids with Gorenstein property. As they wrote the whole classification seems to be `quite difficult'. For other combinatorial objects this question was also intensively studied -- e.g.~for perfect matchings polytopes of grid graphs \cite{BeHaSa09}, for cut polytopes \cite{Oh14}, for symmetric edge polytopes \cite{MaHiNaOhHi11,HiKuMi17}, and for order polytopes \cite{HiMaOhSh15}. For three special classes of matroids a classification of Gorenstein matroids was obtained -- for uniform matroids \cite{DeHi97}, for lattice path matroids \cite{Kn20}, and recently for graphic matroids \cite{HiLaMaMiVo19} (with an extension to multigraphs \cite{Ko19}). We obtain a complete classification of matroids whose variety is Gorenstein.

\subsection{Combinatorial meaning}

Matroid base polytope is the convex hull of the indicator vectors of all bases of the matroid. It is a well-established object of study in matroid theory. Its edges correspond to a single element exchanges between bases \cite{GeGoMaSe87}. Thus, the $1$-skeleton of the matroid base polytope is the basis graph of the matroid -- a graph on all bases of the matroid and edges between two bases differing only by one element. Matroid base polytopes are generalized permutohedra \cite{ArBeDo10}, and possess integer Carath\'{e}odory property \cite{GiRe12}.

Another polytope naturally associated to a matroid is the matroid independence polytope. As the name suggests, it is the convex hull of the indicator vectors of all independent sets.

A lattice polytope $P$ is Gorenstein if there exists a positive integer $\delta$ and a point $v$ in $\delta P$ such that $v$ has the smallest positive distance to every facet of $\delta P$ among all lattice points. In other words, a polytope $P$ is Gorenstein if and only if for some positive integer $\delta$ the polytope $\delta P$ is reflexive. Reflexive polytopes play a prominent role in algebraic combinatorics with important connections to mirror symmetry \cite{Ba94, BaNi08}. A lot of studies concern their classification and properties, cf.~\cite{KrSk00,Sa00}.  

We classify matroids for which the base polytope or the independence polytope is Gorenstein. One of the biggest challenges was that there are plenty of examples of such matroids. In order to capture all of them in Sections \ref{SectionFMatroids} and \ref{SectionGMatroids} we introduce two new families of matroids parameterized by families of subsets with certain intersection properties. We believe that these results, apart from algebraic meaning, are also interesting for their own combinatorial sake. As one of possible applications, our results give a construction of a large class of Gorenstein polytopes on which conjectures involving Gorenstein polytopes may be verified.

\subsection{Our results}

By Definition \ref{DefinitionGorenstein}, Lemmas \ref{LemmaProductIndependence}, \ref{LemmaProductGorenstein}, and Theorem \ref{TheoremFClassification} we obtain a classification of matroids for which the independence polytope is Gorenstein. The class of $F_{\delta}$-matroids (which appear in the classification) is constructed in Section \ref{SectionFMatroids}.

\begin{theorem}\label{TheoremClassificationP}
The independence polytope $P(M)$ of a matroid $M$ is Gorenstein if and only if there exists an integer $\delta\geq 2$ such that matroid $M$ is a direct sum of loops and $(\delta-1)$-blow ups of connected $F_{\delta}$-matroids.
\end{theorem}

By Definition \ref{DefinitionGorenstein}, Lemmas \ref{LemmaProductBase}, \ref{LemmaProductGorenstein}, and Theorems \ref{TheoremGClassification} and \ref{TheoremG2Classification} we obtain a classification of matroids for which the base polytope is Gorenstein. The class of $G'_\delta$-matroids (which appear in the classification) is constructed in Section \ref{SectionGMatroids}.

\begin{theorem}\label{TheoremClassificationB} 
The base polytope $B(M)$ of a matroid $M$ is Gorenstein if and only if there exists an integer $\delta\geq 2$ such that: 
when $\delta>2$, matroid $M$ is a direct sum of loops, coloops, and $G'_\delta$-matroids with some good $(\delta-1)$-ears contracted; 
when $\delta=2$, matroid $M$ is a direct sum of loops, coloops, and $G'_2$-matroids.
\end{theorem}

\section{Polytopes -- definitions and properties}

Throughout the paper by $M$ we denote a matroid, by $E$ its ground set, and by $r$ the rank function $2^E\rightarrow\N$. The matroidal closure of a set $A$ we denote by $cl(A)$. We say that a set $F\subset E$ is a \emph{flat}, or that it is \emph{closed}, when $cl(F)=F$. A set $A$ is called \emph{indecomposable} when it cannot be decomposed into proper subsets $A=A_1\sqcup A_2$ (by $\sqcup$ we denote the disjoint union of sets) such that $r(A)=r(A_1)+r(A_2)$. A set $A$ is called \emph{connected} if every two elements of $A$ belong to some circuit contained in $A$. Recall that a set is connected if and only if it is indecomposable if and only if it is not a direct sum of two or more nontrivial matroids. We will use these notions interchangeably. In particular, a connected matroid is loopless. In Sections \ref{SectionBGorenstein} and \ref{SectionGMatroids} we will intensively use one more notion -- a \emph{flacet} is a flat $G$ such that both: the restriction of the matroid to $G$ and the contraction of $G$ in the matroid are connected. This notion was introduced in \cite{FeSt05}. In \cite{HiLaMaMiVo19} the name \emph{good flat} was used instead. For a general background of matroid theory we refer the reader to \cite{Ox92}. 

\subsection{Independence polytope}

\begin{definition}
	The \emph{independence polytope} of a matroid $M$ on the ground set $E$, denoted by $P(M)$, is the convex hull of points $v_I:=\sum_{e\in I}\chi_e\in\Z^E$ over all independent sets $I$ of the matroid $M$. 
\end{definition}

Indicator vectors of independent sets are vertices of the hypercube $[0,1]^E$, hence every indicator vector $e_I$ of an independent set $I$ is a vertex of $P(M)$. 

It is straightforward to show the following lemma.

\begin{lemma}\label{LemmaProductIndependence}
Suppose a matroid $M$ is the direct sum of matroids $M_1,\dots,M_k$. Then the independence polytope $P(M)$ is the Cartesian product of the independence polytopes $P(M_1),\dots,P(M_k)$. In particular, the independence polytope of a loop is a vertex, so when $e$ is a loop in $M$, then $P(M)$ and $P(M\setminus e)$ are lattice isomorphic. 
\end{lemma} 

In Section \ref{SectionPGorenstein} we will use the following (not as obvious as the above) lemma describing facets of the independence polytope $P(M)$.

\begin{customlemma}{\ref{LemmaFFacets}}
	Let $M$ be a loopless matroid. The matroid independence polytope $P(M)$ is full dimensional in $\R^E$, and the following set of inequalities is minimal defining $P(M)$:
	\begin{itemize}
		\item[(i)] $0\leq x_e$, for every $e\in E$,
		\item[(ii)] $\sum_{e\in F} x_e\leq r(F)$, for every indecomposable flat $F$.
	\end{itemize}
	That is, the intersection of $P(M)$ with each of the supporting hyperplanes of the above half spaces is a facet of $P(M)$.
\end{customlemma}

\subsection{Base polytope}

\begin{definition}
The \emph{base polytope} of a matroid $M$ on the ground set $E$, denoted by $B(M)$, is the convex hull of points $v_B:=\sum_{e\in B}\chi_e\in\Z^E$ over all bases $B$ of the matroid $M$. 
\end{definition}

Since indicator vectors of bases of $M$ are vertices of the hypercube $[0,1]^E$, every indicator vector $v_B$ of a basis $B$ is a vertex of $B(M)$. Clearly, $B(M)=P(M)\cap\{v:\sum_{e\in E}v_e=r(M)\}$, thus $B(M)$ is a face of $P(M)$.

It is straightforward to show the following lemma.

\begin{lemma}\label{LemmaProductBase}
Suppose a matroid $M$ is the direct sum of matroids $M_1,\dots,M_k$. Then the base polytope $B(M)$ of $M$ is the Cartesian product of the base polytopes $B(M_1),\dots,B(M_k)$. In particular, the base polytope of a loop (or a coloop) is a single vertex, so when $e$ is a loop (or a coloop) in $M$, then $B(M)$ and $B(M\setminus e)$ are lattice isomorphic. 
\end{lemma}

In Section \ref{SectionBGorenstein} we will use the following (even less obvious than in the case of the independence polytope) lemma describing facets of the base polytope $B(M)$.

\begin{customlemma}{\ref{LemmaGFacets}}
	Let $M$ be a connected matroid. The matroid base polytope $B(M)$ is full dimensional in an affine hyperplane $L:=\{x\in\R^E:\sum_{e\in E} x_e=r(E)\}$, and the following set of inequalities is minimal defining $B(M)$ in the hyperplane $L$:
	\begin{itemize}
		\item[(i)] $0\leq x_e$, for every $e\in E$ such that $M\setminus e$ is connected,
		\item[(ii)] $\sum_{e\in G} x_e\leq r(G)$, for every proper \emph{flacet} $G$ -- a flat $\emptyset\neq G\subsetneq E$ such that: restriction of $M$ to $G$ and contraction of $G$ in $M$ are connected.
	\end{itemize}
	That is, the intersection of $B(M)$ with each of the supporting hyperplanes of the above half spaces is a facet of $B(M)$.
\end{customlemma}

\subsection{Gorenstein polytope}\label{SubsectionGorensteinPolytope}

Recall that a lattice polytope is \emph{reflexive} if $0$ is the only lattice point in its interior and the dual (polar) polytope is again a lattice polytope. A lattice polytope $P$ has \emph{integer decomposition property} if every vector in $kP\cap\Z^d$ is a sum of $k$ vectors from $P\cap\Z^d$ (cf. \cite{CoHaHibHig14}).

\begin{definition}\label{DefinitionPreGorenstein}
A full-dimensional lattice polytope $P\subset\R^d$ with integer decomposition property is \emph{Gorenstein} if there exists a positive integer $\delta$ and a lattice point $v\in\delta P$ such that $\delta P-v$ is a reflexive polytope.
\end{definition}

The above is a combinatorial definition of the Gorenstein property. It is already in use on its own combinatorial sake -- see e.g. \cite{BaNi08, BeHaSa09, BrRo07, HiHiOh11, NiSc13, Oh14, OhHi06, SaSt18}. However, the notion of a Gorenstein polytope originates from algebra. We will not present here a definition of a Gorenstein algebra. Instead, we list a few equivalent conditions with combinatorial, algebraic, and geometric meaning. Let $P\subset \R^d$ be a full-dimensional lattice polytope with integer decomposition property, and let $\C[P]$ be the semigroup algebra of lattice points in the cone over $P$ -- a standard toric construction \cite{Fu93,CoLiSc11,St95,BrGu09}. The following conditions are equivalent \cite{Ba94,Hi92}:
\begin{enumerate}
	\item $P$ is a Gorenstein polytope,
	\item $\C[P]$ is a Gorenstein algebra, 
	\item the affine variety associated to $\C[P]$ is Gorenstein,
	\item the numerator of the Hilbert series of $\C[P]$ is palindromic,
	\item the canonical divisor of $\Spec \C[P]$ is Cartier.
\end{enumerate}
Moreover, by virtue of the work of Batyrev \cite{Ba94} Gorenstein polytopes play an important role in mirror symmetry and for this reason are very intensively studied, see e.g. \cite[Section 8.3]{CoLiSc11} and references therein.

We will use the following `working' description of the Gorenstein property of a polytope. As we easily show in Proposition \ref{PropositionGorenstein} this coincides with Definition \ref{DefinitionPreGorenstein}.

\begin{definition}\label{DefinitionGorenstein}
Let $\delta$ be a positive integer. A full-dimensional lattice polytope $P\subset \R^d$ with integer decomposition property is \emph{$\delta$-Gorenstein} if there exists a lattice point $v\in\delta P$ such that for every supporting hyperplane of the cone over $P$ its reduced equation $h$ (that is, $h$ such that $h(\Z^d)=\Z$) satisfies $h(v)=1$. 

A non necessarily full-dimensional lattice polytope $P\subset \R^d$ with integer decomposition property is \emph{$\delta$-Gorenstein} if $P$ is $\delta$-Gorenstein in the affine lattice it spans.

A lattice polytope $P$ with integer decomposition property is \emph{Gorenstein} if $P$ is $\delta$-Gorenstein for some positive integer $\delta$.
\end{definition}

\begin{proposition}\label{PropositionGorenstein}
Definitions \ref{DefinitionPreGorenstein} and \ref{DefinitionGorenstein} agree for full-dimensional polytopes.
\end{proposition}

\begin{proof}
	When $P$ is $\delta$-Gorenstein, then vertices of the dual polytope $(\delta P-v)^*$ correspond to reduced equations of supporting hyperplanes of the cone over $P$, and therefore they are integral. Conversely, if $h(v)>1$ for the reduced equation of a supporting hyperplane of the cone over $P$, then the coordinates of the corresponding vertex of the dual polytope $(\delta P-v)^*$ are fractions with denominator $h(v)$, hence they are not integer. 
\end{proof}

\begin{remark}
	The polytopes $P(M)$ and $B(M)$ for a matroid $M$ do not contain interior lattice points, thus they are not $1$-Gorenstein unless they are $0$-dimensional.
\end{remark}

\begin{lemma}\label{LemmaProductGorenstein}
	Let $P_1,P_2$ be two lattice polytopes with integer decomposition property. The product $P_1\times P_2$ is $\delta$-Gorenstein if and only if $P_1$ and $P_2$ are $\delta$-Gorenstein.
\end{lemma}

\begin{proof}
	If both $P_1$ and $P_2$ are $\delta$-Gorenstein with lattice points $v_1$ and $v_2$, as in Definition \ref{DefinitionGorenstein}, then $(v_1,v_2)\in\delta(P_1\times P_2)$ is the point which proves that $P_1\times P_2$ is $\delta$-Gorenstein. Conversely, suppose $P_1\times P_2$ is $\delta$-Gorenstein with a point $(v_1,v_2)$ with lattice distance one from all facets. Then $v_1$ (resp.~$v_2$) is a point which proves that $P_1$ (resp.~$P_2$) is $\delta$-Gorenstein.  
\end{proof}

In particular, when $M$ is a matroid which is the direct sum of matroids $M_i$, then $P(M)$ (resp. $B(M)$) is $\delta$-Gorenstein if and only if for every $i$ the polytope $P(M_i)$ (resp. $B(M_i)$) is $\delta$-Gorenstein.

\section{$F$-matroids}\label{SectionFMatroids}

In this section we define a class of matroids that will play a central role in classification of Gorenstein independence polytopes. We first define a special family of subsets of the ground set -- a family satisfying the properties of indecomposable flats in a matroid whose independence polytope is Gorenstein. Our goal is to achieve a matroid whose set of indecomposable flats coincides with the original family. 

\subsection{$F_{\delta}$-families}\label{SubsectionFFamilies}

Fix an integer $\delta\geq 2$. 

\begin{definition}\label{DefinitionFFamily} 
A family $\mathcal{F}$ of subsets of a finite ground set $E$ is called an \emph{$F_{\delta}$-family} if it satisfies the following conditions:
	\begin{enumerate}
	\item[(1)] for every $e\in E$, $\{e\}\in\mathcal{F}$, 
    \item[(2)] for every $F_i,F_j\in \mathcal{F}$, if $F_i\cap F_j\neq\emptyset$, then $F_i\cap F_j$ and $F_i\cup F_j$ belong to $\mathcal{F}$,
    \item[(3)] for every $F_i\in \mathcal{F}$ its cardinality equals $\frac{\delta f_i-1}{\delta-1}$ for some integer $f_i\equiv 1$ (mod $\delta-1$).
	\end{enumerate}
We call the $F_{\delta}$-family $\mathcal{F}$ \emph{connected} when $E\in\mathcal{F}$.
\end{definition}

Notice that for $F_k=F_i\cup F_j$ and $F_l=F_i\cap F_j$ as in $(2)$, by counting their cardinalities we get an equality $f_k+f_l=f_i+f_j$ between the corresponding numbers.

Moreover, notice that in an $F_{\delta}$-family it cannot happen that $F_k=F_i\sqcup F_j$ (that is, $F_k=F_i\cup F_j$ and $F_i\cap F_j=\emptyset$) for sets $F_k,F_i,F_j$ from the family. Indeed, it is impossible by the congruence of their cardinalities modulo $\delta$. We call it the non-disjoint union of two property.

\smallskip

$F_{\delta}$-families can be easily generated. 

Either globally -- by setting an intersection scheme satisfying condition $(2)$ between sets $F_i$, assigning numbers $f_i$ in a strictly monotone way, and finally filling every set with the right number of points (so that condition $(3)$ holds). Or inductively -- every $F_{\delta}$-family is a union of connected $F_{\delta}$-families on disjoint ground sets. We describe two constructions of a connected $F_\delta$-family from smaller connected $F_\delta$-families. Every connected $F_\delta$-family is achieved in one of these ways.

\begin{definition}
	Suppose $\mathcal{F}_1$ and $\mathcal{F}_2$ are two connected $F_\delta$-families on ground sets $E_1,E_2$ such that $E_1\cap E_2=F\in\mathcal{F}_1\cap\mathcal{F}_2$ and $\{A_1\in\mathcal{F}_1:A_1\subset F\}=\{A_2\in\mathcal{F}_2:A_2\subset F\}$. The \emph{fiber sum} of families $\mathcal{F}_1,\mathcal{F}_2$ is a connected $F_\delta$-family $F$ on the ground set $E_1\cup E_2$ defined by $\mathcal{F}:=\mathcal{F}_1\cup\mathcal{F}_2\cup\{A_1\cup A_2:A_i\in\mathcal{F}_i,A_1\cap F=A_2\cap F\neq\emptyset\}$. We denote it by $\mathcal{F}_1\uplus \mathcal{F}_2$.
\end{definition}

Showing that $\mathcal{F}_1\uplus \mathcal{F}_2$ is indeed an $F_\delta$-family is almost straightforward.  There are several cases to consider, out of which the least obvious is to check that when $(A_1\cup A_2)\cap (A'_1\cup A'_2)\neq\emptyset$ and it is contained in $E_1\setminus E_2$, then the union is in $\mathcal{F}$. But then, $(A_1\cup A'_1)\cap F=(A_1\cap F)\sqcup(A'_1\cap F)$ and all three sets are in $\mathcal{F}_1$, which violates the non-disjoint union of two property. 

\begin{definition}
	Suppose $\mathcal{F}_1,\dots,\mathcal{F}_k$ are $k$ connected $F_\delta$-families on disjoint ground sets $E_1,\dots,E_k$. Suppose $k\equiv1$ (mod $\delta$). The \emph{connection} of $\mathcal{F}_1,\dots,\mathcal{F}_k$ is a connected $F_\delta$-family $\mathcal{F}$ on the ground set $E=E_1\cup\dots\cup E_k$ defined by $\mathcal{F}:=\mathcal{F}_1\cup\dots\cup \mathcal{F}_k\cup\{E\}$. We denote it by $\Conn(\mathcal{F}_1,\dots,\mathcal{F}_k)$.
\end{definition}

\begin{proposition}\label{PropositionFConstruction}
	A connected $F_\delta$-family on $\vert E\vert>1$ is equal to either:
	\begin{itemize}
		\item a fiber sum $\mathcal{F}_1\uplus\mathcal{F}_2$ of connected $F_\delta$-families, or
		\item a connection $\Conn(\mathcal{F}_1,\dots,\mathcal{F}_k)$ of $k\equiv1$ (mod $\delta$) connected $F_\delta$-families.
	\end{itemize}
\end{proposition}

\begin{proof}
	Suppose $\mathcal{F}$ is a connected $F_\delta$-family on $\vert E\vert>1$. Consider maximal (w.r.t. inclusion) proper subsets of $E$ that belong to $\mathcal{F}$. There are two cases:
	\begin{itemize}
		\item there exist two sets with nonempty intersection, or
		\item all such sets are pairwise disjoint.
	\end{itemize}
	
	In the first case, let $F_1,F_2\in\mathcal{F}$ be such that $F_1\cap F_2\neq\emptyset$. Then, by $(2)$ we have that $F_1\cup F_2\in\mathcal{F}$, hence (from maximality of $F_1$ and $F_2$) $F_1\cup F_2=E$. Consider two connected $F_\delta$-families $\mathcal{F}_i=\{F\in\mathcal{F}:F\subset F_i\}$ for $i=1,2$. Then, $\mathcal{F}=\mathcal{F}_1\uplus \mathcal{F}_2$ as a consequence of the non-disjoint union of two property -- every set from $\mathcal{F}$ disjoint from $F_1\cap F_2$ is contained in $F_1$ or in $F_2$.
	
	In the second case, let $F_1,F_2,\dots,F_k$ be all of them. Since $\vert E\vert>1$, $k>1$. Counting cardinalities, we get $k\equiv1$ (mod $\delta$). Consider connected $F_\delta$-families $\mathcal{F}_i=\{F\in\mathcal{F}:F\subset F_i\}$ for $i=1,\dots,k$. Then, $\mathcal{F}=\Conn(\mathcal{F}_1,\dots,\mathcal{F}_k)$.
\end{proof}

\subsection{$F_{\delta}$-matroids}\label{SubsectionFMatroids}

We present a construction of matroids corresponding to $F_{\delta}$-families. 

\begin{definition}\label{DefinitionFMatroid}
Let $\mathcal{F}$ be an $F_{\delta}$-family. Let $\mathfrak{C}_{\mathcal{F}}$ be a family of minimal sets (w.r.t. inclusion) among all $(f_i+1)$-element subsets of $F_i$ over all $F_i\in\mathcal{F}$. The matroid corresponding to $\mathcal{F}$ is a matroid on the ground set $E$ with the set of circuits equal to $\mathfrak{C}_{\mathcal{F}}$. We denote it by $M_{\mathcal{F}}$, and call an \emph{$F_\delta$-matroid}.	
\end{definition}
As we will see in Theorem \ref{TheoremFMatroids}, $M_{\mathcal{F}}$ is indeed a matroid. 

\begin{remark}\label{RemarkF1}
A set $A\subset E$ is independent in $M_{\mathcal{F}}$ if and only if for every $F_i\in\mathcal{F}$ an inequality $\lvert F_i\cap A\rvert\leq f_i$ holds. In particular, the rank of $F_i$ in $M_{\mathcal{F}}$ is less or equal to $f_i$.
\end{remark}

The following is the main theorem of our classification of matroids whose independence polytope is Gorenstein. It shows that conditions from Definition \ref{DefinitionFFamily} are not only  necessary conditions for a family of indecomposable flats (see Subsection \ref{SubsectionPGorensteinCombinatorial}), but also sufficient conditions.

\begin{theorem}\label{TheoremFMatroids}
Suppose $\mathcal{F}$ is an $F_{\delta}$-family for some $\delta\geq 2$. Then, $M_{\mathcal{F}}$ is a matroid in which the set of indecomposable flats (closed and connected sets) is equal to $\mathcal{F}$, and every $F_i\in\mathcal{F}$ has rank equal to $f_i$. In particular, the matroid $M_{\mathcal{F}}$ is connected if and only if $\mathcal{F}$ is a connected $F_{\delta}$-family.
\end{theorem}

\begin{proof}
Firstly, we show that $M_{\mathcal{F}}$ is indeed a matroid, that is we show that the set $\mathfrak{C}_{\mathcal{F}}$ satisfies circuit axioms (see \cite{Ox92}). The first two are clearly satisfied. For the third, let $C_i$ be a circuit, i.e.~a $(f_i+1)$-element subset of $F_i$, and let $C_j$ be another circuit, i.e. a $(f_j+1)$-element subset of $F_j$. Let $x\in C_i\cap C_j$. It is enough to show that there exists a circuit $C_k\subset C_i\cup C_j\setminus x$. Since $x\in F_i\cap F_j$, by $(2)$ we get that $F_l=F_i\cap F_j$ and $F_k=F_i\cup F_j$ belong to $\mathcal{F}$. Now, $\vert C_i\cap C_j\rvert\leq f_l$, as otherwise a $(f_l+1)$-element subset of $F_l$ would be a subset of $C_i$ and $C_j$ which is a contradiction ($C_i$ and $C_j$ would not me minimal). Thus, $\vert C_i\cup C_j\rvert\geq f_i+1+f_j+1-f_l$, and so $\vert C_i\cup C_j\setminus x\rvert\geq f_i+f_j-f_l+1= f_k+1$. It is a subset of $F_k$ of cardinality at least $f_k+1$, so it contains a circuit.
	
\smallskip

Next, we prove that $F_i\in\mathcal{F}$ has rank equal to $f_i$ by induction on the size. If the size of $F_i$ is $1$, then clearly its rank equals $1$ since there are no loops in $M_{\mathcal{F}}$. When the size of $F_i$ is greater than $1$, consider all maximal proper subsets of $F_i$ that belong to $\mathcal{F}$. There are two cases:
	\begin{itemize}
		\item there exist two sets with nonempty intersection, or
		\item all such sets are pairwise disjoint.
	\end{itemize}

In the first case, let $F_j$ be one of them. So, there exists a proper subset $F_k\in\mathcal{F}$ of $F_i$ such that $F_i=F_j\cup F_k$. Let $F_k$ be the minimum set with this property. Let $F_j\cap F_k=F_l$. From maximality of $F_j$ follows that any proper subset $F_x\in\mathcal{F}$ of $F_i$ is either contained in $F_j$, contained in $F_k\setminus F_j$, or it contains $F_k$ (from minimality of $F_k$). Let $I_l$ be a basis of $F_l$ (by ind.~ass.~of size $f_l$), and let $I_j$ (by ind.~ass.~of size $f_j$) and $I_k$ (by ind.~ass.~of size $f_k$) be its extensions to bases of $F_j$ and $F_k$ respectively. Then, the set $I_i=I_j\cup I_k$ of size $f_j+f_k-f_l=f_i$ is independent in $F_i$. Indeed, it is easy to verify that for every set $F_x\in\mathcal{F}$ we have $\lvert F_x\cap I_i \rvert\leq f_x$. If $F_x$ is contained in $F_j$, then it follows from the fact that $I_i\cap F_j=I_j$ is independent. If $F_x$ is contained in $F_k\setminus F_j$, then it follows from the fact that $I_i\cap F_k=I_k$ is independent. Finally, when $F_x$ contains $F_k$, then let $F_y=F_j\cap F_x$. Now, $\lvert F_x\cap I_i \rvert=\lvert F_y\cap I_i \rvert+f_k-f_l\leq f_y+f_k-f_l=f_x$.    
	
In the second case, let $F_1,F_2,\dots,F_k$ be all of them. Clearly, $F_1\sqcup\dots\sqcup F_k=F_i$. Let $I_1,\dots,I_k$ be bases of the corresponding sets. Notice that we have an equality $\frac{\delta f_i-1}{\delta-1}=\vert F_i\vert=\vert F_1\vert+\dots+\vert F_k\vert=\frac{\delta f_1-1}{\delta-1}+\dots+\frac{\delta f_k-1}{\delta-1}$, so $f_i=f_1+\dots+f_k-\frac{(k-1)}{\delta}$. Thus, the only obstructions for the set $I_1\cup\dots\cup I_k$ of size $f_1+\dots+f_k$ to be independent in $F_i$ are the $(f_i+1)$-element circuits of $F_i$. So any $f_i$-element subset of $I_1\cup\dots\cup I_k$ is an independent set in $F_i$, and hence by Remark \ref{RemarkF1} it is a basis. 

\smallskip
	
Now, we argue that sets $F_i\in\mathcal{F}$ are closed. Suppose that $C$ is a circuit in $M_{\mathcal{F}}$ and $\lvert C\setminus F_i\rvert\leq 1$. Since circuits belong to the family $\mathfrak{C}_{\mathcal{F}}$, $C$ is a $(f_j+1)$-element subset of a set $F_j\in\mathcal{F}$. Then $F_k=F_i\cap F_j$ has rank at least $\lvert C\rvert-1=f_j$. But its superset $F_j$ has rank $f_j$. Since both have the same rank, by $(3)$ they have the same cardinality, and hence they coincide. As a consequence, $F_j\subset F_i$. Thus, $C\subset F_i$.
	
\smallskip
	
Notice that the closure of a circuit $C$ in $M_{\mathcal{F}}$ that corresponds to a set $F_i\in\mathcal{F}$ is equal to $F_i$. Indeed, the rank of $cl(C)$ is $\lvert C\rvert-1=f_i$. It is equal to the rank of $F_i$, which is its superset since $F_i$ is closed. Hence, $cl(C)=F_i$. 

\smallskip
	
Finally, let $A$ be a closed set in $M_{\mathcal{F}}$. Consider all maximal (not necessarily proper) subsets of $A$ that belong to $\mathcal{F}$. Denote them by $F_1,\dots,F_k$. They are clearly disjoint (from maximality and $(2)$).  If $k>1$, then the sets $F_1,\dots,F_k$ form a decomposition of $A$ into closed sets in $M_{\mathcal{F}}$. Indeed, otherwise there exists a circuit $C$ corresponding to a set $F\in\mathcal{F}$, which intersects at least two sets $F_i,F_j$ among $F_1,\dots,F_k$. Then $cl(C)=F$ also intersects $F_i,F_j$, and it is contained in $A$ (since $A$ is closed). This contradicts maximality of $F_i$, since $F_i\cup F\cup F_j\in\mathcal{F}$ is a larger subset of $A$. Therefore, if a closed set $A$ is indecomposable, then it belongs to $\mathcal{F}$. 

For the opposite implication, let $F\in\mathcal{F}$. We already know that it is closed. Let $F=F_1\sqcup\dots\sqcup F_k$ be its decomposition into closed and indecomposable subsets (we already know that $F_1,\dots,F_k\in\mathcal{F}$). Then we have two equalities, $f=f_1+\dots+f_k$ and $\vert F\vert=\vert F_1\vert+\dots+\vert F_k\vert$. Together with $(3)$ they give that $k=1$. Thus, $F$ is an indecomposable flat.
\end{proof}

\begin{remark}
Notice that if a connected matroid has the set of indecomposable flats which is an $F_{\delta}$-family it does not imply that it is an $F_{\delta}$-matroid. The uniform matroid $U_{2,5}$ is an example. More generally, the set of indecomposable flats (without their ranks) does not define uniquely a matroid.
\end{remark}

\section{Gorenstein independence polytopes}\label{SectionPGorenstein}

\subsection{Combinatorial reformulation when $P(M)$ is $\delta$-Gorenstein}\label{SubsectionPGorensteinCombinatorial}

Elements $e$ and $f$ of a matroid are called \emph{parallel} when $\{e,f\}$ is a circuit. The relation of being parallel is an equivalence relation in which the equivalence class of $x$ equals $cl(x)$. The following operation allows to enlarge these equivalence classes.

\begin{definition}\label{DefinitionBlowUp}
	A \emph{$k$-blow up of an element $e\in E$ in a matroid $M$} is the matroid $M$ enlarged in the following way: 
	\begin{itemize}
		\item the ground set $E$ is enlarged by new elements $e_2,\dots,e_k$,
		\item the set of bases in enlarged by new bases $(B\setminus e)\cup e_i$ for every $i=2,\dots,k$ and every basis $B$ of $M$ containing $e$.
	\end{itemize}
	Then elements $e,e_2,\dots,e_k$ are parallel elements.
	
	A \emph{$k$-blow up of a matroid} is the $k$-blow up of every element of its ground set.
\end{definition}

\begin{example}
	A $2$-blow up of one element in a graphic matroid of a double edge is a triple edge, while a $2$-blow up of the whole graphic matroid of a double edge is the graphic matroid of four parallel edges.
\end{example}

Notice, that the structure of a matroid in some sense does not change after a $k$-blow up -- as $k$-blow up only makes parallel elements classes $k$ times larger, and being a basis or an independent set depends only on elements' equivalence classes.

\begin{theorem}\label{TranslationP}
	Fix a positive integer $\delta$. Let $M$ be a connected matroid. The matroid independence polytope $P(M)$ is $\delta$-Gorenstein if and only if $M$ is a $(\delta-1)$-blow up of a connected matroid $M'$ which satisfies $(\clubsuit_\delta)$:
	\begin{itemize}
		\item[(1)] $\delta r(F)=(\delta-1)\vert F\vert+1$ for every indecomposable flat $F$ in $M'$.
	\end{itemize}
\end{theorem}

The above theorem is a corollary of Theorem $7.3$ from \cite{HeHi02}. However, for sake of completeness we include a proof. First, we present a description of facets of the matroid independence polytope.

\begin{lemma}[\cite{Sc03}, Theorem $40.5$]\label{LemmaFFacets}
	Let $M$ be a loopless matroid. The matroid independence polytope $P(M)$ is full dimensional in $\R^E$, and the following set of inequalities is minimal defining $P(M)$:
	\begin{itemize}
		\item[(i)] $0\leq x_e$, for every $e\in E$,
		\item[(ii)] $\sum_{e\in F} x_e\leq r(F)$, for every indecomposable flat $F$.
	\end{itemize}
	That is, the intersection of $P(M)$ with each of the supporting hyperplanes of the above half spaces is a facet of $P(M)$.
\end{lemma}

\begin{proof}[Proof of Theorem \ref{TranslationP}]	
	Suppose $(\clubsuit)_{\delta}$ $(1)$ holds for a matroid $M'$ on the ground set $E'$. Let $M$ (on the ground set $E$) be the $(\delta-1)$-blow up of $M'$. Notice that an indecomposable flat $F$ in $M$ is exactly the $(\delta-1)$-blow up of an indecomposable flat $F'$ in $M'$. Further, the rank of $F$ in $M$ coincides with the rank of $F'$ in $M'$. Hence, by $(1)$ we have $\delta r(F)=\delta r(F')=(\delta-1)\vert F'\vert+1=\vert F\vert+1$ for every indecomposable flat $F$ in $M$. Let $v=(1,\dots,1)\in\R^E$. By Lemma \ref{LemmaFFacets} vector $v$ belongs to $\delta P_M$, and further -- reduced equations of facets of $\delta P_M$ of both types $(i),(ii)$ evaluated at $v$ give $1$. Therefore, the polytope $P(M)$ is $\delta$-Gorenstein.
	
	Conversely, suppose the polytope $P(M)$ is $\delta$-Gorenstein. Therefore, there exists a lattice point $v\in\delta P(M)$ such that $\delta P(M)-v$ is a reflexive polytope. Firstly, by Lemma \ref{LemmaFFacets} $(i)$ we have $v_e=1$. By $(ii)$ for every indecomposable flat $F$ in $M$ the equation $\sum_{e\in F} x_e=\delta r(F)$ is the reduced equation of a supporting hyperplane to $\delta P(M)$. Thus, $\vert F\vert+1=\delta r(F)$. For every element $e\in E$ consider its closure $cl(e)$. It is an indecomposable flat of rank $1$. Hence, $\vert cl(e)\vert=\delta-1$. Notice that $cl(e)$ is exactly the set of all elements in $M$ parallel to $e$. Let $M'$ be the matroid $M$ restricted to representatives of all parallel element classes. Then, clearly $M$ is the $(\delta-1)$-blow up of $M'$. Moreover, for every indecomposable flat $F'$ in $M'$ we have $\delta r(F')=\delta r(F)=\vert F\vert+1=(\delta-1)\vert F'\vert+1$. That is, $(\clubsuit)_{\delta}$ $(1)$ holds for $M'$. 
\end{proof}

For the remaining part of this subsection let $\delta\geq 2$ be a fixed integer, and let $M'$ be a fixed matroid on the ground set $E$ satisfying condition $(\clubsuit)_{\delta}$ $(1)$. 

For every $A\subset E$ define
$$c(A):=\delta r(A)-(\delta-1)\vert A\vert.$$ 
It is straightforward to check that the function $c$ has the following properties:
\begin{enumerate}
	\item[(2)] $c(A)=c(cl(A))+(\delta-1)\vert cl(A)\setminus A\vert$,
	\item[(3)] $c(A)\geq c(cl(A))$,
	\item[(4)] $c(A)=c(A_1)+\dots+c(A_k)$ if $A_1,\dots,A_k$ are connected components of $A$,
	\item[(5)] $c(A\cup B)+c(A\cap B)\leq c(A)+c(B)$ -- using submodularity of rank function.
\end{enumerate}
Moreover, using $(\clubsuit)_{\delta}$ $(1)$ it is easy to check that the function $c$ characterizes indecomposable flats:
\begin{enumerate}
	\item[(6)] $c(A)\geq 0$,
	\item[(7)] $c(A)=0$ if and only if $A=\emptyset$,
	\item[(8)] $c(A)=1$ if and only if $A$ is an indecomposable flat.
\end{enumerate}
Indeed, suppose $A\neq\emptyset$. By $(2)$ we get $c(A)=c(cl(A))+(\delta-1)\vert cl(A)\setminus A\vert\geq c(cl(A))$. When the flat $F=cl(A)$ decomposes into $k$ indecomposable flats $F_1,\dots F_k$, then $c(F)=c(F_1)+\dots+c(F_k)=k\geq 1$. Hence, $c(A)\geq 1$. Following these inequalities we get $(8)$.

\begin{lemma}\label{LemmaFFlatsSumInt}
	Indecomposable flats in $M'$ satisfy the following property:
	\begin{enumerate}
		\item[(9)] if $F_1,F_2$ are indecomposable flats and $F_1\cap F_2\neq\emptyset$ then $F_1\cap F_2$ and $F_1\cup F_2$ are also indecomposable flats.
	\end{enumerate}
\end{lemma}

\begin{proof}
	By $(5)$ and $(6)$ we have $0\leq c(F_1\cup F_2)+c(F_1\cap F_2)\leq c(F_1)+c(F_2)=2$. Moreover, $F_1\cap F_2\neq\emptyset$ and $F_1\cup F_2\neq\emptyset$, thus by $(7)$ we get $c(F_1\cup F_2)=c(F_1\cap F_2)=1$, and by $(8)$ sets $F_1\cap F_2,F_1\cup F_2$ are indecomposable flats.
\end{proof}

\subsection{Classification when $P(M)$ is $\delta$-Gorenstein}

The following is our classification of matroids whose independence polytope is $\delta$-Gorenstein. The class of $F_{\delta}$-matroids (which appear in the classification) is constructed in Section \ref{SectionFMatroids}.

\begin{theorem}\label{TheoremFClassification}
	Fix an integer $\delta\geq 2$. The independence polytope $P(M)$ of a connected matroid $M$ is $\delta$-Gorenstein if and only if $M$ is a $(\delta-1)$-blow up of a connected $F_{\delta}$-matroid.
\end{theorem}

\begin{proof}
	Suppose the independence polytope $P(M)$ of a connected matroid $M$ is $\delta$-Gorenstein. Then by Theorem \ref{TranslationP} the matroid $M$ is the $(\delta-1)$-blow up of a connected matroid $M'$ which satisfies $(\clubsuit_\delta)$ $(1)$. Denote the ground set of $M'$ by $E$. Let $\mathcal{F}$ be the set of all indecomposable flats in $M'$. Notice that, due to $(9)$ and $(1)$, $\mathcal{F}$ is a connected $F_{\delta}$-family on the set $E$ (from Definition \ref{DefinitionFFamily}). Now, consider matroids $M'$ and $M_{\mathcal{F}}$ (from Definition \ref{DefinitionFMatroid}). Clearly, both are on the same ground set $E$. Using Theorem \ref{TheoremFMatroids} we get that both matroids have the same rank, the same set of indecomposable flats $\mathcal{F}$, and that the ranks of these indecomposable flats coincide. Therefore, $M'=M_{\mathcal{F}}$. Indeed, by Lemma \ref{LemmaFFacets} the independence polytopes of both matroids are cut by the same set of halfspaces ($0\leq x_e$ over all $e\in E$, and $\sum_{e\in F} x_e\leq r(F)$ over all $F\in\mathcal{F}$). Thus, $P(M')=P(M_{\mathcal{F}})$ and therefore $M'=M_{\mathcal{F}}$ is an $F_{\delta}$-matroid.
	
	Suppose now $M$ is a $(\delta-1)$-blow up of a connected $F_{\delta}$-matroid $M_{\mathcal{F}}$ (for some $F_{\delta}$-family $\mathcal{F}$). By Theorem \ref{TranslationP}, it is enough to show that the connected matroid $M_{\mathcal{F}}$ satisfies condition $(\clubsuit)_{\delta}$ $(1)$. It does -- by Theorem \ref{TheoremFMatroids} $\mathcal{F}$ is the set of indecomposable flats, so $(\clubsuit)_{\delta}$ $(1)$ follows from $(3)$ in Definition \ref{DefinitionFFamily} and again Theorem \ref{TheoremFMatroids}.
\end{proof}

\begin{example}
	Fix an integer $\delta\geq 2$. The independence polytope of the graphic matroid of the $(\delta-1)$-blow up of the $(\delta+1)$-cycle is $\delta$-Gorenstein, see \cite{HiLaMaMiVo19}. By Theorem \ref{TheoremFClassification} the graphic matroid of a $(\delta+1)$-cycle is a connected $F_{\delta}$-matroid. Indeed, the corresponding $F_{\delta}$-family consists of $\delta+1$ singletons $\{e_i\}$ and a set $E=\{e_1,\dots,e_{\delta+1}\}$.
\end{example}

\begin{example}
	Fix an integer $\delta\geq 2$. The independence polytope of the graphic matroid of the $(\delta-1)$-blow up of two $(\delta+1)$-cycles joined by an edge is $\delta$-Gorenstein, see \cite{HiLaMaMiVo19}. By Theorem \ref{TheoremFClassification} the graphic matroid of two $(\delta+1)$-cycles joined by an edge is a connected $F_{\delta}$-matroid. Indeed, the corresponding $F_{\delta}$-family consists of $2\delta+1$ singletons $\{e_i\}$, sets $F_1=\{e_1,\dots,e_{\delta+1}\}$, $F_2=\{e_{\delta+1},\dots,e_{2\delta+1}\}$, and a set $E=\{e_1,\dots,e_{2\delta+1}\}$.
\end{example}

\begin{example}
	The independence polytope of the uniform matroid $U_{4,7}$ is $2$-Gorenstein. Indeed, by Theorem \ref{TheoremFClassification} it is ($1$-blow up of) a connected $F_{2}$-matroid corresponding to the $F_{2}$-family consisting of a $7$-element set and $7$ singletons.
\end{example}

\section{$G$-matroids}\label{SectionGMatroids}

In this section we define a class of matroids that is essential for classification of $\delta$-Gorenstein base polytopes. In analogy to what we did for the independence polytopes, we define a special family of subsets of the ground set -- a family satisfying the properties of flacets in a matroid whose base polytope is $\delta$-Gorenstein. Our goal is to achieve a matroid whose set of flacets coincides with the original family. We distinguish two cases $\delta>2$ and $\delta=2$.

\subsection{$G_\delta$-families for $\delta\geq 2$}\label{SubsectionGFamilies}

Fix an integer $\delta\geq 2$. The results of this subsection are valid both for $\delta>2$ and $\delta=2$, and will be used later in Subsections  \ref{SubsectionGdelta} and \ref{SubsectionG2}.

\begin{definition}\label{DefinitionGFamily}
A family $\mathcal{G}$ of subsets of a finite ground set $E$ is called a \emph{$G_{\delta}$-family} if it satisfies the following conditions:
\begin{enumerate}
	\item[(1)] for every $e\in E$, $\{e\}\in\mathcal{G}$,
	\item[(2)] for every $G_i,G_j\in \mathcal{G}$, if $G_i\cap G_j\neq\emptyset$ and $G_i\cup G_j\neq E$, then $G_i\cap G_j$ and $G_i\cup G_j$ belong to $\mathcal{G}$,
	\item[(3)] the cardinality of $E$ equals $\frac{\delta g_0}{\delta-1}$ for some integer $g_0\equiv 0$ (mod $\delta-1$),
	\item[(4)] for every $G_i\in \mathcal{G}$ its cardinality equals $\frac{\delta g_i-1}{\delta-1}$ for some $g_i\equiv 1$ (mod $\delta-1$).
\end{enumerate}
\end{definition}

$G_{\delta}$-families can be quite easily generated. 

Either globally -- by setting an intersection scheme satisfying condition $(2)$ between sets $G_i$, assigning numbers $g_i$ in a strictly monotone way, and finally filling every set with the right number of points (so that conditions $(3)$ and $(4)$ hold). Or locally -- from connected $F_\delta$-families in some analogy how schemes (resp.~manifolds) are constructed from affine schemes (resp.~open discs).

\begin{definition}
We say that a $G_\delta$-family $\mathcal{G}$ comes from an \emph{atlas} of connected $F_\delta$-families $\mathcal{F}_1,\dots,\mathcal{F}_k$ on $E_1,\dots,E_k$ when the following conditions are satisfied:
	\begin{itemize}
		\item $\mathcal{G}=\bigcup_{i=1}^k \mathcal{F}_i$ -- covering,
		\item $\{A\in\mathcal{F}_i:A\subset E_i\cap E_j\}=\{A\in\mathcal{F}_j:A\subset E_i\cap E_j\}$ -- compatibility.
	\end{itemize}
\end{definition}

\begin{lemma}
	Every $G_\delta$-family $\mathcal{G}$ comes from an atlas of connected $F_\delta$-families.
\end{lemma}

\begin{proof}
When $G\in\mathcal{G}$, then $\mathcal{F}_G:=\{G_i\in\mathcal{G}:G_i\subset G\}$ is clearly a connected $F_\delta$-family -- compare Definitions \ref{DefinitionFFamily} and \ref{DefinitionGFamily}. Let $G_1,\dots,G_k$ be inclusion maximal elements of $\mathcal{G}$. Then $\mathcal{G}$ comes from an atlas of connected $F_\delta$-families $\mathcal{F}_{G_1},\dots,\mathcal{F}_{G_k}$. Indeed, the covering condition follows from the fact that every element of $\mathcal{G}$ is contained in some $G_i$. The compatibility condition is clear, as $F_\delta$-families $\mathcal{F}_{G_i}$ come from a common $G_\delta$-family.
\end{proof}

\begin{remark}
Connected $F_\delta$-families $\mathcal{F}_1,\dots,\mathcal{F}_k$ on $E_1,\dots,E_k$ satisfying the compatibility condition form (an atlas of) a $G_\delta$-family on the ground set $E=\bigcup_{i=1}^k E_i$ by $\mathcal{G}:=\bigcup_{i=1}^k \mathcal{F}_i$ if and only if the following conditions are satisfied:
\begin{itemize}
	\item the cardinality of $E$ equals $\frac{\delta g_0}{\delta-1}$ for some integer $g_0\equiv 0$ (mod $\delta-1$),
	\item if $F_i\in\mathcal{F}_i$, $F_j\in\mathcal{F}_j$, $F_i\cap F_j\neq\emptyset$, and $F_i\cup F_j\neq E$, then $F_i\cup F_j\in\mathcal{F}_s$ for some $s$.
\end{itemize}
\end{remark}

To define a $G_\delta$-matroid we need to build some more structure on a $G_{\delta}$-family.

\smallskip

For a $G_{\delta}$-family $\mathcal{G}$ let $\tilde{\mathcal{G}}$ be the family of nonempty intersections of sets from $\mathcal{G}$.

\smallskip

For $F\in\tilde{\mathcal{G}}$ let $c(F)$ be the least $k$ s.t. $F$ is intersection in $E$ of $k$ sets from $\mathcal{G}$. 

\smallskip

For $F\in\tilde{\mathcal{G}}$ let $r(F)$ be a positive number equal to $\frac{1}{\delta}((\delta-1)\vert F\vert+c(F))$. That is, the following equality $\delta r(F)=(\delta-1)\vert F\vert+c(F)$ holds for $F\in\tilde{\mathcal{G}}$.

\smallskip

We say that an intersection of sets $G_1,\dots,G_k$ is \emph{transversal} when $G_i\cup G_j=E$ for every $i\neq j$. Equivalently, an intersection of $G_1,\dots,G_k\in\mathcal{G}$ is transversal if and only if sets $E\setminus G_1,\dots,E\setminus G_k$ are pairwise disjoint. 

\smallskip

Observe that every $F\in\tilde{\mathcal{G}}$ is a transversal intersection of $c(F)$ sets from $\mathcal{G}$. Indeed, if $F$ is an intersection of the least number of sets from $\mathcal{G}$, then $G_i\cup G_j=E$ for every $i\neq j$. Otherwise, if $G_i\cup G_j\neq E$, we could replace $G_i,G_j$ by $G_i\cap G_j$ which by $(2)$ would belong to $\mathcal{G}$, resulting in a fewer number of sets. 

\smallskip

When $F\in\tilde{\mathcal{G}}$ is a transversal intersection of $k=c(F)$ sets $G_1,\dots,G_k\in\mathcal{G}$, then
\begin{align*}
	\frac{1}{\delta-1}\left(\delta r(F)-k\right)=\vert F\vert= & \;\vert G_1\vert+\dots+\vert G_k\vert-(k-1)\vert E\vert\\
	= & \;\frac{\delta g_1-1}{\delta-1}+\dots+\frac{\delta g_k-1}{\delta-1}-(k-1)\frac{\delta g_0}{\delta-1}\\
	= & \;\frac{1}{\delta-1}\left(\delta(g_1+\dots+g_k-(k-1)g_0)-k\right),
\end{align*}
hence $r(F)=g_1+\dots+g_k-(k-1)g_0$ is an integer.


\smallskip

Of course, $\mathcal{G}\subset\tilde{\mathcal{G}}$. Notice that for $G_i\in\mathcal{G}$ we have $c(G_i)=1$, $r(G_i)=g_i$, and for $E\in\tilde{\mathcal{G}}$ we have $c(E)=0$, $r(E)=g_0$.

\subsection{$G_\delta$-matroids for $\delta\geq 2$}\label{SubsectionGMatroids}

Fix an integer $\delta\geq 2$. The results up to Claim \ref{ClaimGConnected} are valid also for $\delta=2$ and will be used later in Subsection \ref{SubsectionG2}. 

We present a construction of matroids corresponding to $G_{\delta}$-families. The definition is by circuits, but they can be also introduced by the rank function -- using Claim \ref{ClaimGRank}.

\begin{definition}\label{DefinitionGMatroid}
Let $\mathcal{G}$ be a $G_{\delta}$-family. Let $\mathfrak{C}_{\mathcal{G}}$ be a family of minimal sets (w.r.t. inclusion) among all $(r(F_i)+1)$-element subsets of $F_i$ over all $F_i\in\tilde{\mathcal{G}}$. The matroid corresponding to $\mathcal{G}$ is a matroid on the ground set $E$ with the set of circuits equal to $\mathfrak{C}_{\mathcal{G}}$. We denote it by $M_{\mathcal{G}}$, and call a \emph{$G_\delta$-matroid}.	
\end{definition}

\begin{remark}\label{RemarkG}
A set $A\subset E$ is independent in $M_{\mathcal{G}}$ if and only if for every $F_i\in\tilde{\mathcal{G}}$ an inequality $\lvert F_i\cap A\rvert\leq r(F_i)$ holds. In particular, the rank of $F_i$ in $M_{\mathcal{G}}$ is less or equal to $r(F_i)$.
\end{remark}

The following theorem is the cornerstone of our classification of matroids whose base polytope is $\delta$-Gorenstein for $\delta>2$. It shows that conditions from Definition \ref{DefinitionGFamily} are not only a necessary conditions for a family of flacets (see Subsection \ref{SubsectionBGorensteinCombinatorial}), but also sufficient conditions.

\begin{theorem}\label{TheoremGMatroids}
Suppose $\mathcal{G}$ is a $G_{\delta}$-family for some $\delta>2$. Then, $M_{\mathcal{G}}$ is a connected matroid in which the set of flacets (sets that are closed, connected, and their contraction is connected) is equal to $\mathcal{G}$, every $G_i\in\mathcal{G}$ has rank equal to $g_i$, and the ground set $E$ has rank equal to $g_0$. 
\end{theorem}

We prove Theorem \ref{TheoremGMatroids} in a sequence of claims, some of which before proving require some additional lemmas. The whole theorem works for $\delta>2$, but claims and lemmas up to Claim \ref{ClaimGConnected} are valid also for $\delta=2$. Starting from Lemma \ref{LemmaGIncreasing} an additional assumption is made, which holds always when $\delta>2$.

\begin{lemma}\label{LemmaGSubmodularity}
For every $F_i,F_j\in\tilde{\mathcal{G}}$, if $F_i\cap F_j\neq\emptyset$, then $F_i\cap F_j,F_i\cup F_j\in\tilde{\mathcal{G}}$. Moreover, the function $r:\tilde{\mathcal{G}}\rightarrow\N$ is submodular, that is if $F_i,F_j\in\tilde{\mathcal{G}},F_i\cap F_j\neq\emptyset$, then $r(F_i\cap F_j)+r(F_i\cup F_j)\leq r(F_i)+r(F_j)$.
\end{lemma}

\begin{proof}
	Suppose $F_1,F_2\in\tilde{\mathcal{G}}$, that is $E\setminus F_1=E\setminus G^1_1\sqcup\dots\sqcup E\setminus G^1_k$ and $E\setminus F_2=E\setminus G^2_1\sqcup\dots\sqcup E\setminus G^2_l$, where $G^1_i,G^2_j\in\mathcal{G}$, $c(F_1)=k$, and $c(F_2)=l$. 
	
	Suppose that $F_1\cap F_2\neq\emptyset$. Just from the definition of $\tilde{\mathcal{G}}$ we have that the intersection $F_1\cap F_2\in\tilde{\mathcal{G}}$. For the union we have an equality $E\setminus (F_1\cup F_2)=(E\setminus F_1)\cap (E\setminus F_2)=\bigsqcup_{i,j} (E\setminus G^1_i)\cap (E\setminus G^2_j)=\bigsqcup_{i,j:G^1_i\cup G^2_j\neq E} E\setminus (G^1_i\cup G^2_j)$. Since $G^1_i\cap G^2_j\neq\emptyset$ (as $F_1\cap F_2\neq\emptyset$) and $G^1_i\cup G^2_j\neq E$, we have that $G^1_i\cup G^2_j\in\mathcal{G}$. Hence $F_1\cup F_2\in\tilde{\mathcal{G}}$.
	
	The cardinality $\vert F\vert$ is a modular function. From $r(F)=\frac{1}{\delta}((\delta-1)\vert F\vert+c(F))$ we get that the function $r(F)$ is submodular if and only if the function $c(F)$ is. We will show that $c$ is a submodular function. 
	
	Let $H$ be a bipartite graph with two classes of vertices $A=\{a_1,\dots,a_k\}$ and $B=\{b_1,\dots,b_l\}$. An edge joins vertices $a_i$ and $b_j$ if the corresponding sets $E\setminus G^1_i$ and $E\setminus G^2_j$ have nonempty intersection. In other words, $H$ is the intersection graph of sets $E\setminus G^1_1,\dots,E\setminus G^1_k,E\setminus G^2_1,\dots,E\setminus G^2_l$, as sets $E\setminus G^1_i$ over $i$ (and also sets $E\setminus G^2_j$ over $j$) are pairwise disjoint.
	
	Observe that $c(F_1\cap F_2)$ is at most the number of connected components of $H$. Indeed, $E\setminus (F_1\cap F_2)=E\setminus F_1\cup E\setminus F_2=(E\setminus G^1_1\sqcup\dots\sqcup E\setminus G^1_k)\cup(E\setminus G^2_1\sqcup\dots\sqcup E\setminus G^2_l)$. Whenever in this union there are two sets $E\setminus G_i,E\setminus G_j$ with nonempty intersection we replace them by their union $E\setminus (G_i\cap G_j)$ (which is also of a form $E\setminus G$ for some $G\in\mathcal{G}$, as $\emptyset\neq F_1\cap F_2\subset G_i\cap G_j$). This way sets corresponding to vertices of every connected component of $H$ will become their union as edges of $H$ indicate which pairs of sets have nonempty intersection. In general, the union of sets corresponding to vertices of a connected subgraph $H'$ of $H$ is of a form $E\setminus G_{H'}$ for some $G_{H'}\in\mathcal{G}$. 
	
	Our previous bound on $c(F_1\cup F_2)$ was $k\cdot l$ -- the number of sets $(E\setminus G^1_i)\cap (E\setminus G^2_j)$. In order to get a better bound we need to cover all these sets with a fewer number of sets of a form $E\setminus G$ for $G\in\mathcal{G}$, contained in $E\setminus (F_1\cup F_2)=(E\setminus F_1)\cap (E\setminus F_2)$. Notice that when $H'$ is a connected subgraph of $H$, $a_i\notin H'$, and $H'\cap B=\{b_{j_1},\dots,b_{j_s}\}$, then $E\setminus G^1_i\cap E\setminus G_{H'}=(E\setminus G^1_i)\cap (E\setminus G^2_{j_1})\sqcup\dots\sqcup (E\setminus G^1_i)\cap (E\setminus G^2_{j_s})$. Moreover, this set is either empty, or of a form $E\setminus G$ for some $G\in\mathcal{G}$. We say that a \emph{strange pair} $(a_i,H')$ (where $H'$ is a connected subgraph, $a_i\notin H'$, and $H'\cap B=\{b_{j_1},\dots,b_{j_s}\}$) covers $(a_i,b_{j_1}),\dots,(a_i,b_{j_s})$. Our task is to cover all intersections $E\setminus G^1_i\cap E\setminus G^2_j$. But when $a_i$ and $b_j$ are not joined by an edge, this set is empty. Hence, it is enough to cover edges of $H$ by strange pairs of types $(a_i,H')$ and $(H',b_i)$ (analogously). The number of strange pairs in a covering of the edge set of $H$ (which depends purely on the graph $H$) gives an upper bound on $c(F_1\cup F_2)$. 
	
	We show by induction on $n$ that in every graph $G$ on $n$ vertices the sum of the number of connected components and the number of strange pairs in some covering (in the above sense) of all its edges is at most $n$. First, notice that it is enough to consider connected graphs. Second, let $v$ be a vertex which does not disconnect the graph $G$, i.e. $G\setminus v$ is connected (there always exists such a vertex). Then, a strange pair $(v,G\setminus v)$ covers all edges incident to $v$. Moreover, from the inductive assumption the connected graph $G\setminus v$ possesses a covering with at most $n-2$ pairs. Hence, $G$ has a covering with at most $n-1$ pairs and the inductive assertion follows.
	
	Concluding, $c(F_1\cup F_2)+c(F_1\cup F_2)\leq\vert H\vert=k+l=c(F_1)+c(F_2)$. 
\end{proof}

\begin{claim}\label{ClaimGMatroid}
	$M_{\mathcal{G}}$ is a matroid.
\end{claim}

\begin{proof}
Let $C_i$ be a circuit, i.e.~a $(r(F_i)+1)$-element subset of $F_i$, and let $C_j$ be another circuit, i.e. a $(r(F_j)+1)$-element subset of $F_j$. Let $x\in C_i\cap C_j$. From the circuit axioms (see \cite{Ox92}), it is enough to show that there exists a circuit $C_k\subset C_i\cup C_j\setminus x$. Since $x\in F_i\cap F_j$, by Lemma \ref{LemmaGSubmodularity} we get that $F_l=F_i\cap F_j$ and $F_k=F_i\cup F_j$ belong to $\tilde{\mathcal{G}}$. Now, $\vert C_i\cap C_j\rvert\leq r(F_l)$, as otherwise a $(r(F_l)+1)$-element subset of $F_l$ would be a subset of $C_i$ and $C_j$ which is a contradiction ($C_i$ and $C_j$ would not be minimal). Thus, $\vert C_i\cup C_j\rvert\geq r(F_i)+1+r(F_j)+1-r(F_l)$, and so $\vert C_i\cup C_j\setminus x\rvert\geq r(F_i)+r(F_j)-r(F_l)+1\geq r(F_k)+1$ by Lemma \ref{LemmaGSubmodularity}. It is a subset of $F_k$ of cardinality at least $r(F_k)+1$, so it contains a circuit.
\end{proof}

Using the submodular function $r:\tilde{\mathcal{G}}\rightarrow\N$ we define a new function $\tilde{r}$ defined on all subsets of $E$ in a following way -- for $A\subset E$ let
$$\tilde{r}(A)=\min\{r(F_1)+\dots+r(F_k):F_1,\dots,F_k\in\tilde{\mathcal{G}}\text{ and }A\subset F_1\cup\dots\cup F_k\}.$$

\begin{lemma}\label{LemmaGrtilde}
	The function $\tilde{r}:2^{E}\rightarrow\N$ is 
	\begin{enumerate}
		\item proper -- for $X\subset E$, $0\leq\tilde{r}(X)\leq\vert X\vert$, 
		\item weakly increasing -- for $Y\subset X$, $\tilde{r}(Y)\leq\tilde{r}(X)$, and
		\item submodular -- for $X,Y\subset E$, $\tilde{r}(X\cap Y)+\tilde{r}(X\cup Y)\leq\tilde{r}(X)+\tilde{r}(Y)$.
	\end{enumerate}
\end{lemma}

\begin{proof}
    Clearly, $\tilde{r}$ is proper (by \ref{DefinitionGFamily} $(1)$) and weakly increasing. Let $A_1,A_2\subset E$. Suppose that $\tilde{r}(A_1)=r(F^1_1)+\dots+r(F^1_k)$ and $\tilde{r}(A_2)=r(F^2_1)+\dots+r(F^2_l)$, where $F^1_1,\dots,F^1_k,F^2_1,\dots,F^2_l\in\tilde{\mathcal{G}}$, $A_1\subset F^1_1\cup\dots\cup F^1_k$, and $A_2\subset F^2_1\cup\dots\cup F^2_l$. Moreover, we can assume that sets $F^1_1,\dots,F^1_k$ (and also sets $F^2_1,\dots,F^2_l$) are pairwise disjoint (if two sets have nonempty intersection we can replace them by their union). For submodularity of $\tilde{r}$ we need to show that 
	$$\tilde{r}(A_1\cap A_2)+\tilde{r}(A_1\cup A_2)\leq\tilde{r}(A_1)+\tilde{r}(A_2)=r(F^1_1)+\dots+r(F^1_k)+r(F^2_1)+\dots+r(F^2_l).$$
	We begin on the right side of the above formula and apply the following process -- if there are two sets $C_i,C_j$ intersecting properly (their intersection is nonempty and one is not contained in the other) we exchange them into $C_i\cap C_j$ and $C_i\cup C_j$. After each step all sets belong to $\tilde{\mathcal{G}}$, the multiset union of all sets remains the same, and the sum of values of $r$ on these sets weakly decreases, as by Lemma \ref{LemmaGSubmodularity} an inequality $r(C_i\cap C_j)+r(C_i\cup C_j)\leq r(C_i)+r(C_j)$ holds. The process clearly ends, as for e.g. the sum of squares of cardinalities of the sets grows, and at the same time it is bounded. We obtain 
	$$r(C_1)+\dots+r(C_s)\leq r(F^1_1)+\dots+r(F^1_k)+r(F^2_1)+\dots+r(F^2_l),$$
	where $C_1\cup\dots\cup C_s=F^1_1\cup\dots\cup F^1_k\cup F^2_1\cup\dots\cup F^2_l$ as multisets (notice that the multiplicity of every element is either $1$, or $2$), and no $C_i,C_j$ intersect properly. Observe that the last property implies that $C_1,\dots,C_s$ can be split into (without loss of generality) $C_1,\dots,C_t$ and $C_{t+1},\dots,C_s$ such that as sets
	\begin{align*}
		C_1\cup\dots\cup C_t= & \;(F^1_1\cup\dots\cup F^1_k)\cap (F^2_1\cup\dots\cup F^2_l),\\
		C_{t+1}\cup\dots\cup C_s= & \;(F^1_1\cup\dots\cup F^1_k)\cup(F^2_1\cup\dots\cup F^2_l).
	\end{align*}
	Then, we get inequalities 
	\begin{align*}
		\tilde{r}(A_1\cap A_2)\leq & \;\tilde{r}((F^1_1\cup\dots\cup F^1_k)\cap (F^2_1\cup\dots\cup F^2_l))\leq r(C_1)+\dots+r(C_t),\\
		\tilde{r}(A_1\cup A_2)\leq  & \;\tilde{r}((F^1_1\cup\dots\cup F^1_k)\cup(F^2_1\cup\dots\cup F^2_l))\leq r(C_{t+1})+\dots+r(C_s),
	\end{align*}
	and finally
	$$\tilde{r}(A_1\cap A_2)+\tilde{r}(A_1\cup A_2)\leq\tilde{r}(A_1)+\tilde{r}(A_2).$$
\end{proof}

\begin{claim}\label{ClaimGRank}
	The rank in $M_{\mathcal{G}}$ is given by the function $\tilde{r}$.
\end{claim}

\begin{proof}
By \cite[Corollary $1.3.4$]{Ox92} a proper, weakly increasing, and submodular function is the rank function of a matroid. Thus, Lemma \ref{LemmaGrtilde} guarantees that $\tilde{r}$ is the rank function of a matroid. Denote this matroid by $M_{\tilde{r}}$. We will show the set of circuits of $M_{\tilde{r}}$ coincides with the set of circuits of $M_{\mathcal{G}}$, and therefore both matroids are the same.

Suppose $C$ is a circuit of $M_{\mathcal{G}}$. Then, $C$ is a $(r(F_i)+1)$-element subset of $F_i\in\tilde{\mathcal{G}}$. We have $\tilde{r}(C)\leq r(F_i)=\vert C\vert-1<\vert C\vert$ since $C\subset F_i$. Every proper subset $C'\subsetneq C$ is independent in $M_{\mathcal{G}}$, so if $C'\subset F'_1\cup\dots\cup F'_k$, then $\vert C'\vert\leq\vert C'\cap F'_1\vert+\dots\vert C'\cap F'_k\vert\leq r(F'_1)+\dots+r(F'_k)$ by Remark \ref{RemarkG}. Hence, $\vert C'\vert\leq\tilde{r}(C')$. As a consequence, $C$ is a circuit in $M_{\tilde{r}}$.
	
Suppose now that $C$ is a circuit in $M_{\tilde{r}}$. Then, $\tilde{r}(C)<\vert C\vert$. Suppose $\tilde{r}(C)=r(F_1)+\dots+r(F_k)$ for $F_1,\dots,F_k\in\tilde{\mathcal{G}}$ and $C\subset F_1\cup\dots\cup F_k$. For every $i$ we have $C\cap F_i\subset F_i$. Thus either $k=1$, or for every $i$ the set $C\cap F_i\subsetneq C$ is not a circuit in $M_{\tilde{r}}$, and so $r(F_i)\geq\tilde{r}(C\cap F_i)\geq\vert C\cap F_i\vert$. Altogether it gives $r(F_1)+\dots+r(F_k)\geq\vert C\vert$, which is a contradiction. Therefore, $k=1$. Since $C$ is a circuit in $M_{\tilde{r}}$, $r(F_1)=\tilde{r}(C)=\vert C\vert-1$, so $\vert C\vert=r(F_1)+1$. Now, if $C$ was not minimal in the set of $(r(F_1)+1)$-element subsets of $F_1\in\tilde{\mathcal{G}}$, it would also be not minimal in the set $\{C\subset E:\tilde{r}(C)<\vert C\vert\}$. Hence, $C$ is a circuit in $M_{\mathcal{G}}$.
\end{proof}

\begin{remark}
	Notice that it may happen for $F\in\tilde{\mathcal{G}}$ that $\tilde{r}(F)<r(F)$. It may even happen that $\vert F\vert<r(F)$. In particular, $F$ does not have to be the closure of a $(r(F)+1)$-element circuit in $M_{\mathcal{G}}$.
\end{remark}

\begin{claim}\label{ClaimGRank2}
Every $G_i\in\mathcal{G}$ has rank equal to $g_i$, and $E$ has rank equal to $g_0$. Moreover, for every $\emptyset\neq X\subsetneq E$ the rank of $X$ is greater of equal to $\frac{1}{\delta}((\delta-1)\vert X\vert+1)$.
\end{claim}

\begin{proof}
For $\emptyset\neq X\subsetneq E$, by Claim \ref{ClaimGRank} we have that
\begin{align*}
	\tilde{r}(X)= & \;\min\{r(F_1)+\dots+r(F_k):F_i\in\tilde{\mathcal{G}},X\subset \bigcup_iF_i\}\\
                     =& \;\min\left\{\sum_i\frac{1}{\delta}((\delta-1)\vert F_i\vert+c(F_i)):F_i\in\tilde{\mathcal{G}},X\subset\bigcup_iF_i\right\}\\
                \geq& \;\frac{1}{\delta}((\delta-1)\vert X\vert+1).
\end{align*}

For $G_i\in\mathcal{G}$ we get that $\tilde{r}(G_i)\geq\frac{1}{\delta}((\delta-1)\vert G_i\vert+1)=r(G_i)=g_i$. On the other hand, by Remark \ref{RemarkG} $\tilde{r}(G_i)\leq g_i$. Hence, $\tilde{r}(G_i)=g_i$. 

A similar calculation for $E$ gives $\tilde{r}(E)=\frac{1}{\delta}(\delta-1)\vert E\vert=r(E)=g_0$.
\end{proof}

\begin{claim}\label{ClaimGConnected}
	The matroid $M_{\mathcal{G}}$ is connected.
\end{claim}

\begin{proof}
	Let $F_1,\dots,F_k$ be the decomposition of the ground set $E$ into indecomposable flats. Suppose $k>1$. Then, $\emptyset\neq F_i\subsetneq E$ and by Claim \ref{ClaimGRank2}
	\begin{align*}
	\frac{1}{\delta}(\delta-1)\vert E\vert=\tilde{r}(E) & =\tilde{r}(F_1)+\dots+\tilde{r}(F_k)\\
	                                                      & \;\geq\frac{1}{\delta}((\delta-1)\vert F_1\vert+1)+\dots+\frac{1}{\delta}((\delta-1)\vert F_k\vert+1)>\frac{1}{\delta}(\delta-1)\vert E\vert.
	\end{align*}
This is a contradiction. Hence, $k=1$ and therefore the matroid is connected.
\end{proof}

The remaining part of this subsection holds for $G_{\delta}$-families $\mathcal{G}$ satisfying an additional condition:
\begin{enumerate}
	\item[(4$\frac{1}{2}$)] for every $e\in E$, $E\setminus e\notin \mathcal{G}$.
\end{enumerate}

Notice that when $\delta>2$, then $(5)$ already follows from conditions $(3)$ and $(4)$. Indeed, if $G_i=E\setminus e \in\mathcal{G}$, then $(\delta-1)\vert E\vert-(\delta-1)=(\delta-1)\vert E\setminus e\vert=\delta g_i-1$ and $(\delta-1)\vert E\vert=\delta g_0$. Hence, $\delta g_0=\delta g_i+(\delta-1)-1=\delta(g_i+1)-2$, which is possible for integers $g_0,g_i$ only when $\delta=2$.

Suppose now that $\delta>2$, or $\delta=2$ and $E\setminus e\notin \mathcal{G}$ for every $e\in E$. 

\begin{lemma}\label{LemmaGIncreasing}
The function $r$ is strictly increasing, that is if $F_i,F_j\in\tilde{\mathcal{G}},F_i\subsetneq F_j$, then $r(F_i)<r(F_j)$.
\end{lemma}

\begin{proof}
Suppose $F_1,F_2\in\tilde{\mathcal{G}}$, that is $E\setminus F_1=E\setminus G^1_1\sqcup\dots\sqcup E\setminus G^1_k$ and $E\setminus F_2=E\setminus G^2_1\sqcup\dots\sqcup E\setminus G^2_l$, where $G^1_i,G^2_j\in\mathcal{G}$, $c(F_1)=k$, and $c(F_2)=l$. 
	
Suppose $F_1\subsetneq F_2$. To show that $r$ is stricly increasing we need to prove that $(\delta-1)\vert F_1\vert+c(F_1)<(\delta-1)\vert F_2\vert+c(F_2)$ holds. It is equivalent to an inequality $c(F_1)-c(F_2)<(\delta-1)\vert F_2\setminus F_1\vert$. Observe that the number of sets $E\setminus G^1_i$ intersecting nonempty $E\setminus F_2$ is at most $c(F_2)$. Indeed, no more than one of them can intersect nonempty $E\setminus G^2_j$ as otherwise the union of sets which intersect nonempty $E\setminus G^2_j$ together with this set would be of a form $E\setminus G$ for $G\in\mathcal{G}$ allowing a disjoint decomposition of $E\setminus F_1$ into a fewer number of sets. The number of sets $E\setminus G^1_i$ contained in $F_2\setminus F_1$ is less than $\vert F_2\setminus F_1\vert$ as every such set has more than one element (by condition $(5)$). Therefore, the function $r$ is strictly increasing. 
\end{proof}

\begin{claim}\label{ClaimGFlats}
In the matroid $M_{\mathcal{G}}$ there is the following chain of inclusions:
\begin{center}
indecomposable flats (closed and connected sets) $\subset\tilde{\mathcal{G}}\subset$ flats (closed sets).
\end{center}
\end{claim}

\begin{proof}
	First, we argue that a set $F_i\in\tilde{\mathcal{G}}$ is closed. Suppose that $C$ is a circuit in $M_{\mathcal{G}}$ and $\lvert C\setminus F_i\rvert\leq 1$. Since circuits belong to the family $\mathfrak{C}_{\mathcal{G}}$, $C$ is a $(r(F_j)+1)$-element subset of some set $F_j\in\tilde{\mathcal{G}}$. Then, $F_k=F_i\cap F_j\in\tilde{\mathcal{G}}$ has rank at least $\lvert C\rvert-1=r(F_j)$. On the other hand, by Remark \ref{RemarkG} it has rank at most $r(F_k)$. Hence, since by Lemma \ref{LemmaGIncreasing} the function $r$ is strictly increasing, both sets have to coincide. As a consequence, $F_j\subset F_i$. Thus $C\subset F_i$.
	
	Notice that the closure of a circuit $C$ in $M_{\mathcal{G}}$ that corresponds to a set $F_i\in\tilde{\mathcal{G}}$ is equal to $F_i$. Indeed, the rank of $cl(C)$ is $\lvert C\rvert-1=r(F_i)$. It is greater than or equal to the rank of $F_i$, which is its superset since $F_i$ is closed. Hence, $cl(C)=F_i$. 
	
	Finally, let $A$ be a closed set in $M_{\mathcal{G}}$. Consider all maximal subsets of $A$ that belong to $\tilde{\mathcal{G}}$. Denote them by $F_1,\dots,F_k$. They are clearly disjoint (from maximality and Definition \ref{DefinitionGFamily} $(2)$). Moreover, when $k>1$, then sets $F_1,\dots,F_k$ form a decomposition of $A$ into closed sets in $M_{\mathcal{G}}$. That is, $\tilde{r}(A)=\tilde{r}(F_1)+\dots+\tilde{r}(F_k)$. Indeed, otherwise there exists a circuit $C\subset A$ corresponding to a set $F\in\tilde{\mathcal{G}}$, which intersects at least two sets $F_i,F_j$ among $F_1,\dots,F_k$. Then $cl(C)=F\in\tilde{\mathcal{G}}$ also intersects $F_i,F_j$, and it is contained in $A$ (since $A$ is closed). This contradicts maximality of $F_i$, since $F_i\cup F\cup F_j\in\tilde{\mathcal{G}}$ is a larger subset of $A$. Therefore, if a closed set $A$ is indecomposable, then $k=1$ and $A\in\tilde{\mathcal{G}}$. 
\end{proof}

\begin{claim}\label{ClaimGGoodFlats}
The set of flacets in $M_{\mathcal{G}}$ (sets that are closed, connected, and their contraction is connected) is equal to $\mathcal{G}$.
\end{claim}

\begin{proof}
	Suppose $F$ is a flacet in $M_{\mathcal{G}}$. Since $F$ is an indecomposable flat, from Claim \ref{ClaimGFlats} we know that $F\in\tilde{\mathcal{G}}$. Suppose $\tilde{r}(F)=r(F_1)+\dots+r(F_k)$ for $F_i\in\tilde{\mathcal{G}}$, with $F\subset \bigcup_iF_i$. Then, clearly $F=\bigsqcup_iF_i$. Indeed, if $F_i\nsubseteq F$, then it can be replaced by $F_i\cap F$ and $r(F_i\cap F)<r(F_i)$ by Lemma \ref{LemmaGIncreasing}. If $F_i,F_j$ intersect nonempty, then they can be replaced by $F_i\cup F_j$ and $r(F_i\cup F_j)<r(F_i)+r(F_j)$. Now, $\tilde{r}(F)=\tilde{r}(F_1)+\dots+\tilde{r}(F_k)$, and therefore $F$ decomposes into $F_1,\dots,F_k$. Thus, $k=1$. That is, $\tilde{r}(F)=r(F)$.
	
	Suppose now that $c(F)=k$. That is, $F$ is a transversal intersection of $G_1,\dots,G_k\in\mathcal{G}$. Define $F_i=\bigcap_{j\neq i}G_j$. Clearly, $\bigcup_iF_i=E\cup (k-1)F$ as multisets. For every $i$ we have that $c(F_i)=k-1$. Hence, 
	\begin{align*}
		\tilde{r}(F_i)-\tilde{r}(F)\leq r(F_i)-r(F)= & \;\frac{1}{\delta}((\delta-1)\vert F_i\vert+k-1)-\frac{1}{\delta}((\delta-1)\vert F\vert+k)\\
		= & \;\frac{1}{\delta}((\delta-1)\vert F_i\setminus F\vert-1).
	\end{align*}
    Therefore, we get that the sum of ranks of sets $F_i\setminus F$ in $M_{\mathcal{G}}/F$ is less or equal to
	\begin{align*}
	& \;(\tilde{r}(F_1)-\tilde{r}(F))+\dots+(\tilde{r}(F_k)-\tilde{r}(F))\\
	\leq & \;\frac{1}{\delta}((\delta-1)\vert F_1\setminus F\vert-1)+\dots+\frac{1}{\delta}((\delta-1)\vert F_k\setminus F\vert-1)\\
	= & \;\frac{1}{\delta}((\delta-1)\vert E\setminus F\vert-k)=r(E)-r(F)=\tilde{r}(E)-\tilde{r}(F).
    \end{align*}
    Hence, $M_{\mathcal{G}}/F$ decomposes into $F_1,\dots,F_k$. Since $F$ is a flacet, we get that $k=1$. Thus, $F\in\mathcal{G}$.
	
	Other way round, let $G\in\mathcal{G}$. Suppose that $G=F_1\sqcup\dots\sqcup F_k$ is the decomposition of $G$ into indecomposable flats (we already know that $F_1,\dots,F_k\in\mathcal{G}$). Then, by Claim \ref{ClaimGRank2}
	\begin{align*}
		\frac{1}{\delta}((\delta-1)\vert G\vert+1) & \;=\tilde{r}(G)=\tilde{r}(F_1)+\dots+\tilde{r}(F_k)\\
		& \;\geq\frac{1}{\delta}((\delta-1)\vert F_1\vert+1)+\dots+\frac{1}{\delta}((\delta-1)\vert F_k\vert+1)\\
		& \;\geq\frac{1}{\delta}((\delta-1)\vert G\vert+k).
	\end{align*}
    Thus, $k=1$ and so $G$ is an indecomposable flat. 
    
    Let $A_1,\dots,A_k$ be connected components of $M_{\mathcal{G}}/G$. Suppose $k>1$. Then,
    \begin{align*}
    \frac{1}{\delta}((\delta-1)\vert E\setminus G\vert-1)= &\;r(E)-r(G)=\tilde{r}(E)-\tilde{r}(G)\\
    = & \;(\tilde{r}(A_1\cup G)-\tilde{r}(G))+\dots+(\tilde{r}(A_k\cup G)-\tilde{r}(G))\\
    = & \;(\tilde{r}(A_1\cup G)-r(G))+\dots+(\tilde{r}(A_k\cup G)-r(G))\\
    \geq & \;\frac{1}{\delta}(\delta-1)\vert A_1\vert+\dots+\frac{1}{\delta}(\delta-1)\vert A_k\vert\text{ by Claim \ref{ClaimGRank2}}\\
    = & \;\;\frac{1}{\delta}(\delta-1)\vert E\setminus G\vert,
    \end{align*}
    which is a contradiction. Hence, $k=1$ and so $G$ is a flacet.
\end{proof}

\begin{remark}
	A set $F\in\tilde{\mathcal{G}}$ may be decomposable even though $\tilde{r}(F)=r(F)$.
\end{remark}

\subsection{$G'_\delta$-families and $G'_\delta$-matroids for $\delta>2$}\label{SubsectionGdelta}

Fix an integer $\delta >2$. 

\begin{lemma}\label{LemmaDecomposition}
	Suppose $\mathcal{G}$ is a $G_{\delta}$-family. Let $e$ be an element of the ground set $E$ of $M_{\mathcal{G}}$. Then, ${E\setminus e}$ is not connected if and only if $E\setminus e$ decomposes into exactly $\delta-1$ inclusion maximal in $E\setminus e$ sets from $\mathcal{G}$. 
\end{lemma}

\begin{proof}
	Firstly, suppose that $E\setminus e$ decomposes into $G_1,\dots,G_{\delta-1}\in\mathcal{G}$. Then, 
	$$\frac{\delta g_0}{\delta-1}-1=\vert E\setminus e\vert=\vert G_1\vert+\dots+\vert G_{\delta-1}\vert=\frac{\delta g_1-1}{\delta-1}+\dots+\frac{\delta g_{\delta-1}-1}{\delta-1},$$
	so $\tilde{r}((E\setminus e))=\tilde{r}(E)=g_0=g_1+\dots+g_{\delta-1}=\tilde{r}(G_1)+\dots+\tilde{r}(G_{\delta-1})$, and hence $E\setminus e$ decomposes into connected (by Theorem \ref{TheoremGMatroids}) components $G_1,\dots,G_{\delta-1}$. 
	
	Now, suppose that ${E\setminus e}$ decomposes into $k>1$ connected components $C_1,\dots,C_k$. Every component $C_i$ is closed in the matroid $M_{\mathcal{G}}$, as otherwise if $e$ belonged to its closure we would get a decomposition of $M_{\mathcal{G}}$ which is impossible by Theorem \ref{TheoremGMatroids}. By Claim \ref{ClaimGFlats} every component $C_i$ (as it is a closed and connected set) belongs to $\tilde{\mathcal{G}}$, thus $C_i$ is an intersection of some sets $G_i^j$ from $\mathcal{G}$. Clearly, one of these sets $G_i^0$ does not contain $e$, but by Theorem \ref{TheoremGMatroids} $G_i^0$ is connected. Therefore, every component $C_i=G_i^0=:G_i$ belongs to $\mathcal{G}$. We have equalities
	$$\frac{\delta g_0}{\delta-1}-1=\vert E\setminus e\vert=\vert G_1\vert+\dots+\vert G_k\vert=\frac{\delta g_1-1}{\delta-1}+\dots+\frac{\delta g_{k}-1}{\delta-1},$$
	and $g_0=\tilde{r}(E)=\tilde{r}((E\setminus e))=\tilde{r}(G_1)+\dots+\tilde{r}(G_k)=g_1+\dots+g_k$, hence $k=\delta-1$.
\end{proof}

\begin{definition}\label{DefinitionGdeltaFamily}
	A family $\mathcal{G'}$ of subsets of a finite ground set $E$ is called a \emph{$G'_{\delta}$-family} if it is a $G_{\delta}$-family and satisfies the following additional condition:
	\begin{enumerate}
		\item[(5)] for every $e\in E$ the set $E\setminus e$ decomposes into exactly $\delta-1$ inclusion maximal in $E\setminus e$ sets from $\mathcal{G'}$.
	\end{enumerate}
	A \emph{$G'_\delta$-matroid} is the $G_\delta$-matroid of a $G'_{\delta}$-family.
\end{definition}

The following is a direct consequence of the above definition and Lemma \ref{LemmaDecomposition}.

\begin{proposition}\label{PropositionGdeltaMatroid}
	A $G_\delta$-matroid on the ground set $E$ has a property that for every $e\in E$ the set $E\setminus e$ is not connected if and only if it is a $G'_\delta$-matroid.
\end{proposition}

\begin{remark}
$G'_{\delta}$-families can be also quite easily generated, with no need to check condition $(5)$ for $G_{\delta}$-families. Indeed, suppose $\mathcal{G}$ is a $G_\delta$-family on the ground set $E$ such that $e\in E$ is \emph{bad}, that is the set $E\setminus e$ is connected. Let $\mathcal{G'}$ be the family $\mathcal{G}$ in which $e$ is replaced by $(\delta-1)^2$ elements $e_{i,j}$ (i.e. $e$ is replaced in every set from $\mathcal{G}$ to which $e$ belongs) on the ground set $E':=E\setminus\{e\}\cup\{e_{i,j};i,j=1,\dots,\delta-1\}$ and with additional sets $G_i=E'\setminus\{e_{i,1},\dots,e_{i,\delta-1}\}$ for every $i=1,\dots,\delta-1$ added. It is straightforward to check that $\mathcal{G'}$ is a $G_\delta$-family in which every new element $e_{i,j}$ is not bad, and no additional bad elements are created. In particular, after application of the above procedure to every bad element in $\mathcal{G}$ one gets a $G'_{\delta}$-family. In the case of a graphic matroid this replacement procedure corresponds to replacing a bad edge by $\delta-1$ parallel paths each of length $\delta-1$ (in the next section these paths are called $(\delta-1)$-ears).
\end{remark}

\subsection{$(\delta-1)$-ears in $G'_{\delta}$-matroids for $\delta>2$}\label{SubsectionEars}

Fix an integer $\delta >2$. 

In analogy to its outlook in the case of graphic matroids, we call a set of $s$ elements $\{e_1,\dots,e_s\}$ of the ground set of a matroid an \emph{$s$-ear} if every circuit of the matroid contains either none of these elements, or all of them.

A classification of matroids whose base polytope is Gorenstein (see Subsection \ref{SubsectionBClassificationDelta} and Section \ref{SectionSmallRank}) intensively uses $(\delta-1)$-ears in $G'_{\delta}$-matroids. Here we provide a useful characterization of such ears and some basic facts.

\begin{lemma}\label{LemmaEars}
	Let $\mathcal{G}$ be a $G'_\delta$-family on the ground set $E$, and let $A\subset E$. The following conditions are equivalent:
	\begin{enumerate}
		\item $A$ is a $(\delta-1)$-ear in the matroid $M_{\mathcal{G}}$,
		\item every element of $\mathcal{G}$ either contains $A$, is disjoint from $A$, or is a singleton, and $|A|=\delta-1$,
		\item every element of $\mathcal{\tilde G}$ either contains $A$, is disjoint from $A$, or is a singleton, and $|A|=\delta-1$,
		\item $E\setminus A$ is an inclusion maximal element in the $G'_\delta$-family $\mathcal{G}$.
	\end{enumerate}
\end{lemma}

\begin{proof}
	$(1)\Rightarrow (2)$: 
	Let $G\in\mathcal{G}$, then by Theorem \ref{TheoremGMatroids} $G$ is connected as a subset of $M_{\mathcal{G}}$. Suppose $G$ has two distinct elements $a,b$, with $a\in A$. Then, since $G$ is connected, there is a circuit in $G$ containing both $a$ and $b$. From the definition of an ear $G$ contains the whole $A$.
	
	$(2)\Rightarrow (1)$: Taking $a\in A$ and applying condition $(5)$ from Definition  \ref{DefinitionGdeltaFamily} we get that $G:=E\setminus A\in\mathcal{G}$. Now $E=A\sqcup G$. Suppose a circuit $C$ in the matroid $M_{\mathcal{G}}$ intersects $A$. The circuit $C$ comes from some set $F\in \mathcal{\tilde{G}}$, and by condition $(2)$ we get that $A\subset F$. From the existence of the circuit $C$ for $F\in \mathcal{\tilde{G}}$ we get that
	$$\tilde{r}(F)=|C|-1=r(F)=\frac{1}{\delta}((\delta-1)|F|+c(F)).$$
	For a set $F\cap G\in\mathcal{\tilde{G}}$ we have
	$$|C\cap G|\leq\tilde{r}(F\cap G)\leq r(F\cap G)=\frac{1}{\delta}((\delta-1)|F\cap G|+c(F\cap G))$$
	and since $|F\cap G|=|F|-|F\cap A|\leq |F|-(\delta-1)$ and $c(F\cap G)\leq c(F)+1$ we get
	\begin{align*}
		|C\cap G|\leq & \;\frac{1}{\delta}((\delta-1)(|F|-(\delta-1))+c(F)+1)\\
		= & \;\frac{1}{\delta}((\delta-1)|F|+c(F))-\delta+2=\tilde{r}(F)-\delta+2=|C|-\delta+1.
	\end{align*}
	Therefore, $|C\cap A|\geq\delta-1$, and thus $A\subset C$, meaning that $A$ is an ear.
	
	$(2)\Leftrightarrow(3)$ is clear.
	
	$(1)\Rightarrow (4)$: 
	Let $A$ be a $(\delta-1)$-ear. For any $a\in A$ the set $E\setminus a$ decomposes into $\delta-2$ singletons of $A\setminus a$ and the set $E\setminus A$. By condition $(5)$ from Definition \ref{DefinitionGdeltaFamily} we get that $E\setminus A\in\mathcal{G}$ and that it is inclusion maximal in $E\setminus a$. Since $a\in A$ can be chosen arbitrarily, $E\setminus A$ is inclusion maximal in $\mathcal{G}$.
	
	$(4)\Rightarrow (2)$: 
	Let $E\setminus A$ be inclusion maximal in $\mathcal{G}$. The family $\mathcal{G}$ on the complement, that is the family $\{G\in\mathcal{G}:G\subset A\}$, is an fact a $F_\delta$-family as no two sets cover the whole $E$. Therefore, $A$ splits into inclusion maximal sets $G_1,\dots,G_k$ in $\{G\in\mathcal{G}:G\subset A\}$. Thus $E=E\setminus A\sqcup G_1\sqcup\dots\sqcup G_k$.
	
	We claim that every $G_i$ is a singleton. Suppose that this is not the case (without loss of generality) for $G_1$. Take $e\in G_1$. The inclusion maximal elements of $\mathcal{G}$ contained in $E\setminus e$ must be $A, G^1,\dots,G^s,G_2,\dots,G_k$, where $G^1,\dots,G^s$ form a decomposition of $G_1\setminus e$. As $|G^i|\equiv 1$ (mod $\delta$) and $|G_1|\equiv 1$ (mod $\delta$), if $s>0$ then $s\geq\delta$. But this contradicts condition $(5)$ from Definition \ref{DefinitionGdeltaFamily}. Thus, the sets $G_i$ are singletons. Now again by condition $(5)$ from Definition \ref{DefinitionGdeltaFamily}, it follows that $k=\delta-1$ and therefore $A$ is a $(\delta-1)$-ear.
\end{proof}

We call a $G'_\delta$-family \emph{trivial} if it consists only of singletons. Then, by Lemma \ref{LemmaDecomposition} there are exactly $\delta$ singletons, and the corresponding matroid is the graphic matroid of a $\delta$-cycle (see Example \ref{ExampleDeltaCycle}). 

\begin{proposition}\label{PropositionEars}
	In a nontrivial $G'_\delta$-family the following is true:
    \begin{enumerate}
    \item $(\delta-1)$-ears are pairwise disjoint,
	\item every $(\delta-1)$-ear is contained in some set from the family,
	\item no set of the family contains all $(\delta-1)$-ears,
	\item there are at least two $(\delta-1)$-ears.
	\end{enumerate}
\end{proposition}

\begin{proof}
	Let $\mathcal{G}$ be a nontrivial $G'_\delta$-family on $E$, and let $A\subset E$ be a $(\delta-1)$-ear.
	
	$(1)$:
	Suppose $A'$ is a $(\delta-1)$-ear different from $A$ and $a\in A\cap A'$. Then $E\setminus A\in\mathcal{G}$ intersects properly $A'$ which by Lemma \ref{LemmaEars} $(2)$ is a contradiction.
	
	$(2)$:
	Take $e\in E\setminus A$. This is possible as $|E\setminus A|\equiv 1$ (mod $\delta$). The family $\mathcal{G}$ is nontrivial, so $E\setminus e$ is larger than $A$. By condition $(5)$ from Definition \ref{DefinitionGdeltaFamily}, there must be a set $G\in\mathcal{G}$ which is not a singleton, that intersects $A$. By Lemma \ref{LemmaEars} $(2)$ this set contains $A$. 

	$(3)$:
	If such a set existed, we could choose one that is inclusion maximal in $\mathcal{G}$. By Lemma \ref{LemmaEars}, its complement would be a $(\delta-1)$-ear, which is a contradiction.
	
	$(4)$:
	For cardinality reason $E\notin\mathcal{G}$. Hence, as the whole $E$ is covered by sets from the family $\mathcal{G}$, there are at least two inclusion maximal elements in $\mathcal{G}$. By Lemma \ref{LemmaEars} $(4)$ these correspond to $(\delta-1)$-ears.
\end{proof}

\subsection{$G'_2$-families and $G'_2$-matroids}\label{SubsectionG2}

\begin{definition}\label{DefinitionG2Family}
A family $\mathcal{G'}$ of subsets of a finite ground set $E$ is called a \emph{$G'_{2}$-family} if it satisfies the following conditions:
	\begin{enumerate}
		\item[(1')] for every $e\in E$, $\{e\}\in\mathcal{G'}$ and $E\setminus e\notin \mathcal{G'}$,
		\item[(2')] for every $G_i,G_j\in \mathcal{G'}$, if $G_i\cap G_j\neq\emptyset$ and $G_i\cup G_j\neq E$, then $G_i\cap G_j\in\mathcal{G'}$, and if $G_i\cap G_j\neq\emptyset$, $G_i\cup G_j\neq E$, and $G_i\cup G_j\neq E\setminus e$, then $G_i\cup G_j\in\mathcal{G'}$,
		\item[(3)] the cardinality of $E$ equals $2g_0$ for some integer $g_0$,
		\item[(4)] for every $G_i\in \mathcal{G'}$ its cardinality equals $2g_i-1$ for some integer $g_i$.
	\end{enumerate}
\end{definition}
Notice that conditions $(3)$ and $(4)$ are the same as in Definition \ref{DefinitionGFamily} for $\delta=2$, while in conditions $(1')$ and $(2')$ there is a slight difference -- sets $E\setminus e$ for $e\in E$ are treated differently.
 
\smallskip 
 
For a $G'_{2}$-family $\mathcal{G'}$ let $\tilde{\mathcal{G'}}$ be the family of nonempty intersections of sets from $\mathcal{G'}$, and let $c'(F)$, $r'(F)$ for $F\in\tilde{\mathcal{G'}}$ be defined as in Subsection \ref{SubsectionGFamilies}.

\smallskip

We present a construction of matroids corresponding to a $G'_{2}$-family $\mathcal{G'}$. The definition is by circuits, but they can be also introduced by the rank function -- by Proposition \ref{PropositionG2Equal} and using Claim \ref{ClaimGRank}.

\begin{definition}\label{DefinitionG2Matroid}
	Let $\mathcal{G'}$ be a $G'_{2}$-family. Let $\mathfrak{C}_{\mathcal{G'}}$ be a family of minimal sets (w.r.t. inclusion) among all $(r'(F_i)+1)$-element subsets of $F_i$ over all $F_i\in\tilde{\mathcal{G'}}$. The matroid corresponding to $\mathcal{G'}$ is a matroid on the ground set $E$ with the set of circuits equal to $\mathfrak{C}_{\mathcal{G'}}$. We denote it by $M_{\mathcal{G'}}$, and call a \emph{$G'_2$-matroid}.	
\end{definition}

\begin{proposition}\label{PropositionG2Equal}
	Suppose $\mathcal{G'}$ is a $G'_{2}$-family. Then, $\mathcal{G}=\mathcal{G'}\cup\{E\setminus e:e\in E\}$ is a $G_{2}$-family (from Definition \ref{DefinitionGFamily}), and the matroid $M_{\mathcal{G'}}$ coincides with the matroid $M_{\mathcal{G}}$ (from Definition \ref{DefinitionGMatroid}).
\end{proposition}

\begin{proof}
It is straightforward that $\mathcal{G}=\mathcal{G'}\cup\{E\setminus e:e\in E\}$ is a $G_{2}$-family, that is $\mathcal{G}$ satisfies conditions $(1)-(4)$ from Definition \ref{DefinitionGFamily}. Let $\tilde{\mathcal{G}}$ be a family of nonempty intersections of sets from $\mathcal{G}$, and let $c(F)$, $r(F)$ for $F\in\tilde{\mathcal{G}}$ be defined as in Subsection \ref{SubsectionGMatroids}.
	
Clearly, $\tilde{\mathcal{G'}}\subset\tilde{\mathcal{G}}$. Notice that on $F\in\tilde{\mathcal{G'}}$ functions $c'$ and $c$, and therefore also functions $r'$ and $r$, coincide. Indeed, if $E\setminus F=E\setminus G_1\cup\dots\cup E\setminus G_k$ with $G_i\in\mathcal{G'}$, then allowing to use also sets $E\setminus (E\setminus e)=\{e\}$ cannot make $k$ smaller. Moreover, if for $F\in\tilde{\mathcal{G}}$ we have $E\setminus F=E\setminus G_1\sqcup\dots\sqcup E\setminus G_l\sqcup\{e_{l+1}\}\sqcup\dots\sqcup\{e_k\}$ with $G_i\in\mathcal{G'}$ and $k=c(F)$, then $F\subset F'=G_1\cap\dots\cap G_l\in\tilde{\mathcal{G'}}$ and $c(F')=l$. Therefore, $$r(F')=\frac{1}{2}(\vert F'\vert+c(F'))=\frac{1}{2}(\vert F\vert+(k-l)+l)=\frac{1}{2}(\vert F\vert+k)=\frac{1}{2}(\vert F\vert+c(F))=r(F).$$

Suppose $C$ is a $(r(F_i)+1)$-element subset of $F_i$ for $F_i\in\tilde{\mathcal{G}}$. Then, by the above, $C$ is also a $(r'(F'_i)+1)=(r(F_i)+1)$-element subset of $F'_i\supset F_i$ for $F'_i\in\tilde{\mathcal{G'}}$. The opposite follows from $\tilde{\mathcal{G'}}\subset\tilde{\mathcal{G}}$. Hence, $\mathfrak{C}_{\mathcal{G'}}=\mathfrak{C}_{\mathcal{G}}$, and finally $M_{\mathcal{G'}}=M_{\mathcal{G}}$.
\end{proof}

Similarly to the case of $\delta>2$, the following theorem is the main ingredient of our classification of matroids whose base polytope is $2$-Gorenstein. It shows that conditions from Definition \ref{DefinitionG2Family} are not only the necessary conditions for a family of flacets (see Subsection \ref{SubsectionBGorensteinCombinatorial}), but also sufficient conditions.

\begin{theorem}\label{TheoremG2Matroids}
	Suppose $\mathcal{G'}$ is a $G'_{2}$-family. Then, $M_{\mathcal{G'}}$ is a connected matroid in which the set of flacets (sets that are closed, connected, and their contraction is connected) is equal to $\mathcal{G'}$, every $G_i\in\mathcal{G'}$ has rank equal to $g_i$, and the ground set $E$ has rank equal to $g_0$. 
\end{theorem}

The proof of Theorem \ref{TheoremG2Matroids} follows the lines of the proof of Theorem \ref{TheoremGMatroids}. All lemmas, remarks, and claims up to Claim \ref{ClaimGConnected} are valid also for the matroid $M_{\mathcal{G}}$ corresponding to the $G_{2}$-family $\mathcal{G}=\mathcal{G'}\cup\{E\setminus e:e\in E\}$. By Proposition \ref{PropositionG2Equal} $M_{\mathcal{G}}$ is equal to $M_{\mathcal{G'}}$. Lemma \ref{LemmaGIncreasing} and Claims \ref{ClaimGFlats}, \ref{ClaimGGoodFlats} hold for $M_{\mathcal{G'}}$ as for every $e\in E$, $E\setminus e\notin \mathcal{G'}$ by $(1')$ (and proofs of these lemma and claims do not use the condition about the union of sets in $\mathcal{G'}$).

\section{Gorenstein base polytopes}\label{SectionBGorenstein}

\subsection{Combinatorial reformulation when $B(M)$ is $\delta$-Gorenstein}\label{SubsectionBGorensteinCombinatorial}

\begin{theorem}\label{TranslationB}
	Fix a positive integer $\delta$. Let $M$ be a connected matroid. The matroid base polytope $B(M)$ is $\delta$-Gorenstein if and only if $M$ satisfies $(\spadesuit)_{\delta}$:
	\begin{enumerate}
		\item[(0)] $M$ possesses a weight $w:E\rightarrow\{1,\delta-1\}$ satisfying both conditions: \newline
		$w(e) = \left\{ \begin{array}{rcl} 
		1 & \mbox{if} & M\setminus e\text{ is connected,} \\
		\delta-1 & \mbox{if} & M/e\text{ is connected,}
		\end{array}\right.$
		\newline
		in particular if for some $e$ both $M\setminus e,M/e$ are connected, then $\delta=2$,
		\item[(1)] $w(E)=\delta r(E)$,
		\item[(2)] $w(G)+1=\delta r(G)$ for every \emph{flacet} $G\subset E$, i.e.~a flat such that both: restriction of $M$ to $G$ and contraction of $G$ in $M$ are connected.
	\end{enumerate}
\end{theorem}

Before we proceed to the proof of Theorem \ref{TranslationB}, recall a description of facets of the matroid base polytope.

\begin{lemma}[\cite{Ki10,FeSt05}]\label{LemmaGFacets}
	Let $M$ be a connected matroid. The matroid base polytope $B(M)$ is full dimensional in the hyperplane $L:=\{x\in\R^E:\sum_{e\in E} x_e=r(E)\}$, and the following set of inequalities is minimal defining $B(M)$ in the hyperplane $L$:
	\begin{itemize}
	\item[(i)] $0\leq x_e$, for every $e\in E$ such that $M\setminus e$ is connected,
	\item[(ii)] $\sum_{e\in G} x_e\leq r(G)$, for every proper \emph{flacet} $G$ -- a flat $\emptyset\neq G\subsetneq E$ such that: restriction of $M$ to $G$ and contraction of $G$ in $M$ are connected.
\end{itemize}
	That is, the intersection of $B(M)$ with each of the supporting hyperplanes of the above half spaces is a facet of $B(M)$.
\end{lemma}

\begin{proof}[Proof of Theorem \ref{TranslationB}]	
	Suppose $(\spadesuit)_{\delta}$ $(0)$ holds, i.e. the weight function $w$ exists. Let $v$ be a lattice point given by $v_e=w(e)$. The affine hyperplane $L$ and inequalities $(i)$, $(ii)$ define $B(M)$, thus by multiplying by $\delta$ their constants we get an affine hyperplane $L_{\delta}$ and inequalities $(i)_{\delta}$, $(ii)_{\delta}$ defining the dilated polytope $\delta B(M)$. Now, conditions $(\spadesuit)_{\delta}$ give that $v\in\delta B(M)$. Further, we claim that both inequalities $(i)_{\delta}$, $(ii)_{\delta}$ provide reduced equations of the facets of $\delta B(M)$. Condition $(\spadesuit)_{\delta}$ $(0)$ gives that $(i)_{\delta}$ evaluated at $v$ is equal to $1$, and condition $(\spadesuit)_{\delta}$ $(2)$ gives that $(ii)_{\delta}$ evaluated at $v$ is equal to $1$. Therefore, the polytope $B(M)$ is $\delta$-Gorenstein.
	
	Conversely, suppose the polytope $B(M)$ is $\delta$-Gorenstein. Therefore, there exists a lattice point $v\in\delta B(M)$ such that $\delta B(M)-v$ is a reflexive polytope. We will show that a weight function defined by $w(e)=v_e$ satisfies conditions $(\spadesuit)_{\delta}$. First, $v\in\delta B(M)$, so $\sum_{e\in E}v_e=\delta r(E)$ ($(1)$ holds). By $(ii)_{\delta}$ for every flacet $G$ the equation $\sum_{e\in G} x_e=\delta r(G)$ is a reduced equation of a supporting hyperplane to $\delta B(M)$. Thus, $\sum_{e\in G} v_e+1=\delta r(G)$ ($(2)$ holds). Now, if $M\setminus e$ is connected, then by $(i)_{\delta}$ the equation $0=x_e$ is a reduced equation of a supporting hyperplane to $\delta B(M)$. Hence $v_e=1$ (the first part of $(0)$ holds). If $M/e$ is connected, then $\{e\}$ is a flacet. The corresponding supporting hyperplane to $\delta B(M)$ is $x_e\leq \delta r(e)=\delta$. Thus, $v_e+1=\delta$ (the second part of $(0)$ holds).
\end{proof}

The weight function $w:E\rightarrow\{1,\delta-1\}$ from Theorem \ref{TranslationB} is already defined by $(\spadesuit)_{\delta}$ $(0)$ for every $e\in E$. 

\begin{lemma}[\cite{Ox92}, Theorem $4.3.1$]\label{LemmaConnectedRestrictionOrContraction}
	Suppose $M$ is a connected matroid. Then for every element $e$ of the ground set, its deletion $M\setminus e$ is connected or its contraction $M/e$ is connected.
\end{lemma}

\begin{lemma}\label{LemmaDecGood}
	Let $M$ be a connected matroid, let $\emptyset\neq F\subsetneq E$ be an indecomposable flat, and let $F_1,\dots,F_k$ be connected components of $M/F$. Then, $E\setminus F_1,\dots,E\setminus F_k$ are flacets. 
\end{lemma}

\begin{proof}
	From the closure properties $cl_M(E\setminus F_i)=cl_{M/F}((E\setminus F_i)\setminus F)\cup F$. The set $(E\setminus F_i)\setminus F$ is a union of connected components in $M/F$, so it is closed. Hence, $E\setminus F_i$ is closed in $M$, i.e. it is a flat.
	
	The set $F_i\cup F$ is connected in $M$. Indeed, suppose contrary $F_i\cup F=C\sqcup D$. The set $F$ is indecomposable, so without loss of generality we have $F\subset C$. Now, $F_i=C\setminus F\sqcup D$ in $M/F$, so it is a decomposition of a connected component (which is not possible), unless $C=F$. But if $C=F$, then $M$ decomposes into $D$ and $E\setminus D$ contradicting connectivity of $M$. Now, $E\setminus F_i$ is a union of connected sets $F_j\cup F$ (for $j\neq i$) with nonempty intersection $F$, hence it is also connected.  
	
	The matroid $M/(E\setminus F_i)$ is isomorphic to $M/F\vert_{F_i}$, hence it is connected.
	
	Concluding, $E\setminus F_1,\dots,E\setminus F_k$ are flacets in $M$.
\end{proof}

\begin{lemma}\label{LemmaDual}
	The matroid base polytope of $M$ is $\delta$-Gorenstein if and only if the matroid base polytope of $M^*$ is $\delta$-Gorenstein (with weights $1$ and $\delta-1$ reversed).
\end{lemma}

\begin{proof}
	Recall that $B(M^*)=(1,\dots,1)-B(M)$. Hence, $B(M)$ and $B(M^*)$ are isomorphic as lattice polytopes. Moreover, $M\setminus e$ is connected if and only if $(M\setminus e)^*=M^*/e$ is connected.
\end{proof}

For the remaining part of this subsection let $\delta\geq 2$ be a fixed integer, and let $M$ be a fixed matroid satisfying conditions $(\spadesuit)_{\delta}$ $(0),(1),(2)$. 

For every $A\subset E$ define
$$c(A):=\delta r(A)-w(A).$$ 
It is straightforward to check that the function $c$ has the following properties:
\begin{enumerate}
	\item[(3)] $c(A)=c(cl(A))+w(cl(A)\setminus A)$,
	\item[(4)] $c(A)\geq c(cl(A))$,
	\item[(5)] $c(A)=c(A_1)+\dots+c(A_k)$ if $A_1,\dots,A_k$ are connected components of $A$,
	\item[(6)] $c(A)=\frac{1}{k-1}(c(A_1)+\dots+c(A_k))$ if $A_1\setminus A,\dots,A_k\setminus A$ are connected components of $M/A$ (we assume here that $A_i$ contains $A$) -- using $(1)$,
	\item[(7)] $c(A)=c(E\setminus A_1)+\dots+c(E\setminus A_k)$ if $A_1,\dots,A_k$ are connected components of $M/A$ (we assume here that every $A_i$ is disjoint from $A$) -- using $(1)$,
	\item[(8)] $c(A\cup B)+c(A\cap B)\leq c(A)+c(B)$ -- using submodularity of rank function.
\end{enumerate}
Moreover, using $(\spadesuit)_{\delta}$ we prove that the function $c$ characterizes flacets:
\begin{enumerate}
	\item[(9)] $c(A)\geq 0$,
	\item[(10)] $c(A)=0$ if and only if $A=E$ or $A=\emptyset$,
	\item[(11)] $c(A)=1$ if and only if $A$ is a flacet or $A=E\setminus e$ and $w(e)=1$.
\end{enumerate}
Indeed, suppose $A\neq \emptyset,E$. By $(3)$ $c(A)=c(cl(A))+w(cl(A)\setminus A)$. If $F=cl(A)\neq E$, then let $F_1,\dots F_k$ be its decomposition. By $(5)$ we have $c(F)=c(F_1)+\dots+c(F_k)$. Now, let $F_i^1,\dots,F_i^{l_i}$ be connected components of $M/F_i$. Then, by Lemma \ref{LemmaDecGood} every $E\setminus F_i^j$ is a flacet, and therefore by $(2)$ we have $c(E\setminus F_i^j)=1$. By $(7)$ $c(F_i)=c(E\setminus F_i^1)+\dots+c(E\setminus F_i^{l_i})=l_i\geq 1$. Hence, $c(A)\geq 1$, and following these inequalities we get also $(11)$.

\begin{lemma}\label{LemmaGFlatsSumInt}
	Flacets satisfy the following properties:
	\begin{enumerate}
		\item[(12)] if $G_1,G_2$ are flacets and $\left\{\begin{array}{l} 
		G_1\cap G_2\neq\emptyset\\
		G_1\cup G_2\neq E \end{array}\right.$ then $G_1\cap G_2$ is a flacet,
		\item[(13)] if $G_1,G_2$ are flacets and $\left\{\begin{array}{l} 
		G_1\cap G_2\neq\emptyset\\
		G_1\cup G_2\neq E\\
		G_1\cup G_2\neq E\setminus e\text{ or }w(e)\neq 1\end{array}\right.$ then $G_1\cup G_2$ is a flacet.
	\end{enumerate}
\end{lemma}

\begin{proof}
	From $(8)$, $(9)$, and $(11)$ we get that 
	$$0\leq c(G_1\cup G_2)+c(G_1\cap G_2)\leq c(G_1)+c(G_2)=2.$$
	If $G_1\cap G_2\neq\emptyset$ and $G_1\cup G_2\neq E$, then from $(10)$ we know that $1\leq c(G_1\cap G_2)$ and $1\leq c(G_1\cup G_2)$. Hence $c(G_1\cap G_2)=c(G_1\cup G_2)=1$. Thus from $(11)$ the set $G_1\cap G_2$ is a flacet. Also from $(11)$ if $G_1\cup G_2\neq E\setminus e$ or $w(e)\neq 1$ (for every $e$), then $G_1\cup G_2$ is a flacet. 
\end{proof}

\subsection{Classification when $B(M)$ is $\delta$-Gorenstein for $\delta>2$}\label{SubsectionBClassificationDelta}

Recall from Subsection \ref{SubsectionEars} that we call a set of $s$ elements $\{e_1,\dots,e_s\}$ of the ground set of a matroid an \emph{$s$-ear} if every circuit of the matroid contains either none of these elements, or all of them.

\begin{definition}\label{DefinitionAddEar}
	A \emph{matroid $M$ with an element $e\in E$ replaced by an $s$-ear} is the matroid $M$ modified in the following way: 
	\begin{itemize}
		\item the ground set $E$ is enlarged by new elements $e_2,\dots,e_s$,
		\item in the set of circuits, every circuit containing $e$ is replaced by a circuit containing $e=e_1,e_2,\dots,e_s$.
	\end{itemize}
\end{definition}

Clearly, $\{e_1,\dots,e_s\}$ is an $s$-ear in the above matroid. 

This operation is a composition of $s-1$ operations known as \emph{series extension}, see \cite{Ox92}. We can also define a matroid with an element replaced by an $s$-ear as a composition of better known operations. Let $M'^*$ denote the matroid obtained by adding $(s-1)$ elements $e_2,\dots,e_s$ parallel to $e=e_1$ in the matroid $M^*$. Then, $M'=(M'^*)^*$ is the matroid $M$ with $e$ replaced by an $s$-ear. Notice that when the matroid $M$ is graphic, then it is just the replacement of an edge $e$ by a path of $s$ edges $e=e_1,\dots,e_s$ -- in \cite{HiLaMaMiVo19} we called it an $s$-ear because it looks like an ear. 

\begin{definition}\label{DefinitionContractEar}
	Let $M$ be a matroid with an $s$-ear $\{e_1,\dots,e_s\}$. A \emph{matroid $M$ with contracted $s$-ear $\{e_1,\dots,e_s\}$} is the matroid $M/\{e_2,\dots,e_s\}$. 
\end{definition}

In other words, a matroid with a contracted $s$-ear is the matroid with contracted all but one elements of that ear. Notice that operations from Definitions \ref{DefinitionAddEar} and \ref{DefinitionContractEar} are inverse to each other.

\begin{definition}\label{DefinitionBadEar}
An ear in a connected matroid is called \emph{bad} if after contracting all elements of the ear the matroid remains connected. An ear is called \emph{good} if it is not bad, i.e.~when contraction of all elements of the ear disconnects the matroid.
\end{definition}

\begin{proposition}\label{Proposition1ToDelta}
	Fix an integer $\delta>2$. Let $M$ be a connected matroid satisfying conditions $(\spadesuit)_{\delta}$. Suppose $w_M(e)=1$ for an element $e$ of the ground set. Let $M'$ be the matroid obtained from $M$ by replacing $e$ with a $(\delta-1)$-ear $\{e=e_1,\dots,e_{\delta-1}\}$. Then, $M'$ is connected, satisfies conditions $(\spadesuit)_{\delta}$, the weight of every $e_i$ equals $\delta-1$ and the obtained ear is good. 
	
	Conversely, let $M'$ be a connected matroid satisfying conditions $(\spadesuit)_{\delta}$ and fix a good $(\delta-1)$-ear. Then, the matroid $M$ obtained from $M'$ by contracting this $(\delta-1)$-ear is also connected and satisfies conditions $(\spadesuit)_{\delta}$.
\end{proposition}

\begin{proof}
	Suppose $w_M(e)=1$ for some $e$ (so $M\setminus e$ is connected and $M/e$ is not connected). Then $w_{M^*}(e)=\delta-1$ (so $M^*\setminus e$ is not connected and $M^*/e$ is connected). By Lemma \ref{LemmaDual} the dual matroid $M^*$ also satisfies $(\spadesuit)_{\delta}$. Let $M'^*$ denote the matroid $M^*$ with added $(\delta-2)$ parallel elements $e_2,\dots,e_{\delta-1}$ to $e=e_1$ -- so that the set consisting of $e$ and all its parallel elements have cardinality $\delta-1$. Now, since $\delta-1>1$, $M'^*\setminus e_i$ is connected and $M'^*/e_i$ is not connected. Hence, we set $w_{M'^*}(e_i)=1$. It is easy to check that flacets $G'$ in $M'^*$ correspond to flacets $G$ in $M^*$ via the rules that if $e\in G$ then $e_1,\dots,e_{\delta-1}\in G'$, and if $e_i\in G'$ then $e\in G$. Taking into account weights, $M'^*$ satisfies $(\spadesuit)_{\delta}$. Let $M'$ be the dual matroid to $M'^*$, which by Lemma \ref{LemmaDual} satisfies $(\spadesuit)_{\delta}$. Clearly, $M'$ is equal to the matroid $M$ with an element $e$ replaced by a $(\delta-1)$-ear $\{e=e_1,\dots,e_{\delta-1}\}$, and $w_{M'}(e_i)=\delta-1$ for every $e_i$. The opposite implication goes analogously.
\end{proof}

\begin{proposition}\label{PropositionAllToDelta}
	Fix an integer $\delta>2$. Let $M$ be a connected matroid satisfying conditions $(\spadesuit)_{\delta}$. Then, there exists a connected matroid $M'$ satisfying conditions $(\spadesuit)_{\delta}$ with all weights equal to $\delta-1$, such that $M$ is equal to $M'$ with  some good $(\delta-1)$-ears contracted. 
\end{proposition}

\begin{proof}
	Applying Proposition \ref{Proposition1ToDelta} to all elements of weight $1$ in $M$ we get a matroid $M'$ satisfying conditions $(\spadesuit)_{\delta}$ with all weights equal to $\delta-1$. The opposite procedure (to get back from $M'$ to $M$) is by contractions from Definition \ref{DefinitionContractEar}. 
\end{proof}

The following is our classification of matroids whose base polytope is $\delta$-Go\-ren\-stein, for $\delta>2$. The class of $G_\delta$-matroids (which appear in the classification) is constructed in Section \ref{SectionGMatroids}.

\begin{theorem}\label{TheoremGClassification}
	Fix an integer $\delta>2$. The base polytope $B(M)$ of a connected matroid $M$ is $\delta$-Gorenstein if and only if $M$ is a $G'_{\delta}$-matroid with some good $(\delta-1)$-ears contracted.
\end{theorem}

\begin{proof}
	By Theorem \ref{TranslationB} the base polytope $B(M)$ of a connected matroid $M$ is $\delta$-Gorenstein if and only if $M$ satisfies conditions $(\spadesuit)_{\delta}$. 
	
	Suppose a connected matroid $M$ satisfies conditions $(\spadesuit)_{\delta}$. Due to Proposition \ref{PropositionAllToDelta} there exists a connected matroid $M'$ satisfying conditions $(\spadesuit)_{\delta}$ with all weights equal to $\delta-1$, such that $M$ is equal to $M'$ with some good $(\delta-1)$-ears contracted. It is enough to prove that $M'$ is a $G'_\delta$-matroid. Denote the ground set of $M'$ by $E$. Notice that the set $\mathcal{G}$ of flacets in $M'$ is a $G_{\delta}$-family on the set $E$ (from Definition \ref{DefinitionGFamily}). Indeed, since all weights in $M'$ are equal to $\delta-1$, for every $e\in E$ the contraction $M'/e$ is connected and therefore $\{e\}$ is a flacet in $M'$ -- condition $(1)$ holds. Condition $(2)$ follows from Lemma \ref{LemmaGFlatsSumInt} and the fact that all weights are equal to $\delta-1>1$. Condition $(3)$ follows from $(\spadesuit)_{\delta}$ $(1)$ and weight $\equiv\delta-1$. And, finally condition $(4)$ follows from $(\spadesuit)_{\delta}$ $(2)$ and weight $\equiv\delta-1$. Now, consider matroids $M'$ and $M_{\mathcal{G}}$ (from Definition \ref{DefinitionGMatroid}). Clearly, both are on the same ground set $E$. Using Theorem \ref{TheoremGMatroids} we get that both matroids have the same rank, the same set of flacets $\mathcal{G}$, and that the ranks of these flacets coincide. Therefore, $M'=M_{\mathcal{G}}$. Indeed, base polytopes of both matroids are contained in the same affine hyperplane (defined by the rank), and by Lemma \ref{LemmaGFacets} both are cut by the same set of halfspaces ($0\leq x_e$ over all $e\in E$, and $\sum_{e\in G} x_e\leq r(G)$ over all $G\in\mathcal{G}$). Thus, $B(M')=B(M_{\mathcal{G}})$ and therefore $M'=M_{\mathcal{G}}$. Now, by the first part of $(\spadesuit)_{\delta}$ $(0)$ and the fact that weight $\equiv\delta-1$, we get that for every $e\in E$ the set $E\setminus e$ is not connected. Thus, by Proposition \ref{PropositionGdeltaMatroid} $M'$ is a $G'_\delta$-matroid.
	
	Suppose now $M$ is a $G'_\delta$-matroid $M_{\mathcal{G'}}$ (for some $G'_{\delta}$-family $\mathcal{G'}$) with some good $(\delta-1)$-ears contracted. By the second part of Proposition \ref{Proposition1ToDelta}, it is enough to show that the connected matroid $M_{\mathcal{G'}}$ satisfies conditions $(\spadesuit)_{\delta}$. It does -- by Theorem \ref{TheoremGMatroids} $\mathcal{G'}$ is the set of flacets. In particular, for every element $e$ of $M_{\mathcal{G'}}$ by Definition \ref{DefinitionGFamily} $(1)$ the set $\{e\}$ is a flacet, hence $M_{\mathcal{G'}}/e$ is connected and all weights are equal to $\delta-1$. Moreover, by Proposition \ref{PropositionGdeltaMatroid} $M_{\mathcal{G'}}\setminus e$ is not connected. Thus, $(0)$ holds. Now, equations $(\spadesuit)_{\delta}$ $(1)$ and $(2)$ follow from conditions $(3)$ and $(4)$ of the $G'_{\delta}$-family.
\end{proof}

\begin{example}\label{ExampleDeltaCycle}
	Fix an integer $\delta>2$. The base polytope of the graphic matroid of the $\delta$-cycle is $\delta$-Gorenstein, see \cite{HiLaMaMiVo19}. By Theorem \ref{TheoremGClassification} it is a $G'_{\delta}$-matroid with some good $(\delta-1)$-ears contracted. Indeed, it is a $G'_{\delta}$-matroid corresponding to a $G'_{\delta}$-family consisting of $\delta$ singletons $\{e_i\}$ on a set $E=\{e_1,\dots,e_{\delta}\}$.
\end{example}

\begin{example}\label{exm:graphicjoin}
	Fix an integer $\delta>2$. The base polytope of the graphic matroid of the $\delta-1$ disjoint $\delta$-cycles joined by an edge is $\delta$-Gorenstein, see \cite{HiLaMaMiVo19}. By Theorem \ref{TheoremGClassification} it is a $G_{\delta}$-matroid with some good $(\delta-1)$-ears contracted. Indeed, it is a $G'_{\delta}$-matroid corresponding to a $G'_{\delta}$-family on a set $E=\{e_{1,1},\dots,e_{\delta,\delta-1}\}$ consisting of $\delta(\delta-1)$ singletons $\{e_{i,j}\}$, and $\delta$ sets $G_i=E\setminus\{e_{i,1},\dots,e_{i,\delta-1}\}$ with contracted one good $(\delta-1)$-ear $\{e_{1,1},\dots,e_{1,\delta-1}\}$.
\end{example}

\begin{example}
 Let $\delta=3$ and let $E$ be a set with $12$ elements divided into $6$ pairs. It is straightforward to check that the $6$ element family consisting of complements of every pair is a $G_3'$-family. In the corresponding $G_3'$-matroid the circuits are precisely subsets that are unions of exactly $3$ pairs of elements. This is not a graphic matroid. Indeed, by contracting one element in every pair, we obtain the uniform matroid $U^2_6$ which is clearly nongraphic. Notice that all $(\delta-1)$-ears of this matroid, that is the $6$ pairs of elements, are bad ears. This example may be easily modified to give for any $\delta>2$ a nongraphic $G'_{\delta}$-matroid. 
\end{example}

\begin{remark}
The distinction between good and bad $(\delta-1)$-ears may be obtained on the level of $G'_\delta$-families as follows. Fix a set of elements $A$ corresponding to a $(\delta-1)$-ear in a  $G'_\delta$-matroid $M$. We know that the base polytope of the dual matroid $M^*$ is Gorenstein, however the weights of edges in $M^*$ are all equal to $1$ by Lemma \ref{LemmaDual}. Consider the matroid $(M^*)'$ by replacing every element of $M^*\setminus A$ by a $(\delta-1)$-ear. 
Notice that $A$ is good if and only if removing $A$ from $M^*$ disconnects the matroid if and only if removing $A$ from $(M^*)'$ disconnects the matroid. Let $(M^*)''$ be obtained from $(M^*)'$ by removing all elements from $A$ apart from one $a\in A$. Then $A$ is good if and only if the base polytope of $(M^*)''$ is Gorenstein if and only if removing element $a$ disconnects the matroid. By Lemma \ref{LemmaDecomposition} this happens if and only if the maximal elements of the corresponding $G'_\delta$-family for $(M^*)''$ partition $(M^*)''\setminus a$ into $(\delta-1)$ sets. 
\end{remark}

\subsection{Classification when $B(M)$ is $2$-Gorenstein}\label{SubsectionBClassification2}

The following is our classification of matroids whose base polytope is $2$-Go\-ren\-stein. The class of $G'_2$-matroids (which appear in the classification) is constructed in Section \ref{SectionGMatroids}.

\begin{theorem}\label{TheoremG2Classification}
	The base polytope $B(M)$ of a connected matroid $M$ is $2$-Gorenstein if and only if $M$ is a $G'_2$-matroid.
\end{theorem}

\begin{proof}
	By Theorem \ref{TranslationB} the base polytope $B(M)$ of a connected matroid $M$ is $2$-Gorenstein if and only if $M$ satisfies conditions $(\spadesuit)_{2}$. 
	
	Suppose a connected matroid $M$ satisfies conditions $(\spadesuit)_{2}$. Denote the ground set of $M$ by $E$. Let $\mathcal{G}$ be the set of all flacets in $M$. Notice that the set $\mathcal{G'}=\mathcal{G}\cup\{\{e\}:e\in E\}$ is a $G'_{2}$-family on the set $E$ (from Definition \ref{DefinitionG2Family}). Indeed, since $M$ is connected no set $E\setminus e$ is a flat in $M$ -- condition $(1')$ holds. Condition $(2')$ follows from Lemma \ref{LemmaGFlatsSumInt}. Conditions $(3)$ and $(4)$ follow from $(\spadesuit)_{2}$ $(1)$ and $(2)$ and the fact that all weights are equal to $1$. Now, consider matroids $M$ and $M_{\mathcal{G'}}$ (from Definition \ref{DefinitionG2Matroid}). Clearly, both are on the same ground set $E$. Using Theorem \ref{TheoremG2Matroids} we get that both matroids have the same rank, the same set of flacets larger than singletons, namely $\mathcal{G}\setminus\{\{e\}:e\in E\}=\mathcal{G'}\setminus\{\{e\}:e\in E\}$, and that the ranks of these flacets coincide. Therefore, $M=M_{\mathcal{G}}$. Indeed, base polytopes of both matroids are contained in the same affine hyperplane (defined by the rank), and by Lemma \ref{LemmaGFacets} both are cut by the same set of halfspaces ($0\leq x_e$ and $x_e\leq 1$ over all $e\in E$, and $\sum_{e\in G} x_e\leq r(G)$ over all $G\in\mathcal{G}\setminus\{\{e\}:e\in E\}$). Recall that when a flacet is a singleton, then the corresponding supporting hyperplane is $x_e\leq r(e)=1$. Thus, $B(M)=B(M_{\mathcal{G'}})$ and therefore $M=M_{\mathcal{G'}}$.
	
	Suppose now $M$ is a $G'_2$-matroid $M_{\mathcal{G}}$ (for some $G'_2$-family $\mathcal{G}$). By Theorem \ref{TheoremG2Matroids} the set of flacets of the matroid $M_{\mathcal{G}}$ is equal to $\mathcal{G}$. It is easy to verify that the connected matroid $M_{\mathcal{G}}$ satisfies conditions $(\spadesuit)_{2}$. Indeed, equations $(\spadesuit)_{2}$ $(1)$ and $(2)$ follow from conditions $(3)$ and $(4)$ of the $G'_{\delta}$-family $\mathcal{G}$.
\end{proof}

\begin{example}
	The base polytope of the graphic matroid of the clique $K_4$ is $2$-Gorenstein, see \cite{HiLaMaMiVo19}. By Theorem \ref{TheoremG2Classification} it is a $G'_2$-matroid. Indeed, it is a $G'_2$-matroid corresponding to a $G'_2$-family on the set $\{e_{1,2},e_{1,3},e_{1,4},e_{2,3},e_{2,4},e_{3,4}\}$ of edges of $K_4$, consisting of six singletons, and four $3$-element sets corresponding to triangles in the clique $K_4$.
\end{example}

\begin{example}
	The base polytope of the graphic matroid of two cliques $K_4$ joined by an edge which is removed is $2$-Gorenstein, see \cite{HiLaMaMiVo19}. By Theorem \ref{TheoremG2Classification} it is a $G'_2$-matroid. Indeed, it is a $G'_2$-matroid corresponding to a $G'_2$-family on the set $\{e_{1,2},e_{1,3},e_{1,4},e_{2,3},e_{2,4},e_{5,6},e_{3,5},e_{4,5},e_{3,6},e_{4,6}\}$, consisting of ten singletons, four $3$-element sets $\{e_{1,2},e_{1,3},e_{2,3}\}$, $\{e_{1,2},e_{1,4},e_{2,4}\}$, $\{e_{3,5},e_{3,6},e_{5,6}\}$, $\{e_{4,5},e_{4,6},e_{5,6}\}$, two $5$-element sets $\{e_{1,2},e_{1,3},e_{1,4},e_{2,3},e_{2,4}\}$, $\{e_{5,6},e_{3,5},e_{4,5},e_{3,6},e_{4,6}\}$, and four $7$-element sets $\{e_{1,2},e_{1,3},e_{1,4},e_{2,3},e_{2,4},e_{3,5},e_{4,5}\}$, $\{e_{1,2},e_{1,3},e_{1,4},e_{2,3},e_{2,4},e_{3,6},e_{4,6}\}$, $\{e_{5,6},e_{3,5},e_{4,5},e_{3,6},e_{4,6},e_{1,3},e_{1,4}\}$, $\{e_{5,6},e_{3,5},e_{4,5},e_{3,6},e_{4,6},e_{2,3},e_{2,4}\}$.
\end{example}

\begin{example} 
The family of singletons on the $4$-element ground set constitutes a $G'_2$-family. The corresponding $G'_2$-matroid is the uniform matroid $U^2_4$. Thus, the matroid is nonbinary, in particular nongraphic, but its base polytope is Gorenstein.  
\end{example}

\section{Gorenstein base polytopes of small rank}\label{SectionSmallRank}

For a nonempty finite set $E$ there is a unique connected matroid on $E$ of rank $1$. The base polytope of this matroid is Gorenstein. As an application of our results, below we provide a full, finite classification of connected matroids of rank $2$ and $3$ whose base polytope is Gorenstein. It would be interesting to know if similar results may be achieved for any fixed rank.

By Lemma \ref{LemmaDual} matroids of rank $r$ whose base polytope is $\delta$-Gorenstein are in bijection with matroids of corank $r$ whose base polytope is $\delta$-Gorenstein. Notice that contraction of a good $(\delta-1)$-ear decreases the rank by $\delta-2$ and it also decreases the cardinality of the ground set by $\delta-2$. In particular, it does not change the corank. Therefore, using our main Theorem \ref{TheoremClassificationB}, in order to classify connected matroids of rank $r$ whose base polytope is Gorenstein it is enough to:

\begin{enumerate}
\item classify all $G'_\delta$-matroids of corank $r$,
\item take the duals of good $(\delta-1)$-ear contractions of matroids from point $(1)$.
\end{enumerate}

Suppose a $G'_\delta$-family on the ground set $E$ gives rise to a matroid of corank $r$. By condition $(3)$ from Definition \ref{DefinitionGFamily} of a $G_\delta$-family we have that $|E|=r\delta$. Recall from Proposition \ref{PropositionEars} that in a nontrivial $G'_\delta$-family $(\delta-1)$-ears are disjoint. In particular, if the number of $(\delta-1)$-ears is $k$, then $k(\delta-1)\leq r\delta$. For $\delta>2$ this gives $k\leq \frac{3}{2}r$. 

For a $G'_\delta$-family we define its \emph{head} to be the set of elements of the ground set that do not belong to any $(\delta-1)$-ear. Equivalently, in a $G'_\delta$-matroid the head is the intersection of all inclusion maximal flacets.

\subsection{Gorenstein base polytopes of rank two}

\begin{theorem}
There are exactly $5$ connected matroids of rank $2$ whose base polytope is Gorenstein. One is $U^2_4$, which is nongraphic. The other $4$ are graphic, they are dual to good $2$-ear contractions in the $G'_3$-matroid corresponding to the $G'_3$-family consisting of singletons and sets $\{a_1,a_2,b_1,b_2\},\{a_1,a_2,c_1,c_2\},\{b_1,b_2,c_1,c_2\}$ on the ground set $E=\{a_1,a_2,b_1,b_2,c_1,c_2\}$.
\end{theorem}

\begin{proof}
By the discussion at the beginning of the section, we only have to classify $G'_\delta$-families with corank $2$. 

First consider $\delta=2$. Then $|E|=4$. As sets from a $G'_\delta$-family must be odd and of cardinality smaller than $3$, these must be exactly singletons. This gives us $U^2_4$.

Now consider $\delta>2$. We know that $|E|=2\delta$ and by Proposition \ref{PropositionEars} that there are at least two disjoint $(\delta-1)$-ears. Their complements are $(\delta+1)$-element
sets from the $G'_\delta$-family. Due to condition $(2)$ from Lemma \ref{LemmaEars} for these $(\delta-1)$-ears, unless $\delta=3$, there are no other $(\delta+1)$-element sets in the $G'_\delta$-family. Then, the $G'_\delta$-family consists of singletons and two $(\delta+1)$-element sets which gives a contradiction with condition $(5)$ from Definition \ref{DefinitionGdeltaFamily} for two elements from the head of the family. Thus, $\delta=3$ (so $|E|=6$) and there is a third $2$-ear for which there is only one possibility. This already determines the matroid from the statement.
\end{proof}

\subsection{Gorenstein base polytopes of rank three}

\begin{theorem}
	There are exactly $17$ connected matroids of rank $3$ whose base polytope is Gorenstein. For $6$ of them it is $2$-Gorenstein, one of which is graphic. For the next $6$ matroids the base polytope is $3$-Gorenstein -- they are duals to good $2$-ears contractions in a single $G'_3$-matroid on the ground set of size $9$. The base polytope of the last $5$ matroids is $4$-Gorenstein. They are duals to good $3$-ears contractions in a single $G'_4$-matroid on the ground set of size $12$.
\end{theorem}

\begin{proof}
First we show that there are exactly $6$ connected matroids with corank $3$ whose base polytope is $2$-Gorenstein. We have $|E|=6$. The $G'_2$-family is determined by the set of subsets of cardinality $3$. The condition is that they cannot intersect in more than one point. We get one matroid $U^3_6$, where there are no such subsets, one matroid when there is one such subset, two for two subsets (depending if the two subsets are disjoint or not), one for three such subsets and one for four such subsets. Only the last one is graphic, it corresponds to $K_4$. Indeed, if such a matroid is graphic, the graph must have $4$ vertices, as the rank is $3$, and $6$ edges.

Next we classify connected $G'_\delta$-matroids with corank $3$ for $\delta>2$. Denote the number of $(\delta-1)$-ears in the $G'_\delta$-family by $k$. We know that $2\leq k\leq \frac{3}{2}\cdot 3=4,5$. Thus, there are three cases:

$k=4$: As $k(\delta-1)\leq 3\delta$ we obtain $\delta\leq 4$. First, suppose $\delta=4$. Then, $|E|=12$ and $E$ consists of four $3$-ears. In particular, the head is empty. Notice that sets $G_i$ in the $G'_\delta$-family are singletons, or have cardinality $(\delta+1)=5$ or $(2\delta +1)=9$. But, as they must consist of $3$-ears in this case, we can exclude $|G_i|=5$. The sets of cardinality $9$ are determined by ears: there are four of them. We thus obtain one matroid of this type. Now, suppose $\delta=3$. Then $|E|=9$, we have four $2$-ears $A_1,A_2,A_3,A_4$ and the head, which is a singleton $\{x\}$. As maximal elements of the $G'_3$-family are determined (complements of $A_i$), we only have to determine the structure of subsets from our family of cardinality $\delta+1=4$. By condition $(2)$ from Lemma \ref{LemmaEars} and parity reasons, none of those may contain $x$. By condition (5) from Definition \ref{DefinitionGdeltaFamily} applied to $x$, the sets $A_1,A_2,A_3,A_4$ must break into two blocks, say $A_1\cup A_2$ and $A_3\cup A_4$ that are in our $G'_3$-family. This already gives us a correct matroid. We cannot have more subsets of cadinality $4$, as then they would intersect one of the two existing subsets in an ear, which cannot be in our family. 

$k=3$: We have three ears $A_1,A_2,A_3$ and the head $\{a,b,c\}$. Consider the inclusion maximum set $S_i$ in the $G'_\delta$-family that contains $A_i$ and does not contain $a$. It must have cardinality $(\delta+1)$. If both $S_1$ and $S_2$ contain just one ear, they must intersect in $\{b,c\}$. This is not possible, as they must be disjoint. Hence, we may assume that $S_1$ contains $A_1$ and $A_2$. This is only possible if $2(\delta-1)\leq \delta+1$, i.e.~$\delta=3$. But, then $S_3=A_3\cup\{b,c\}$. Replacing $a$ with $b$, we get the set $S_3'=A_3\cup\{a,c\}$. But, the sum $S_3\cup S_3'$ cannot be in the family which is a contradiction.

$k=2$: In this case we have two ears $A_1,A_2$ and the head of cardinality $\delta+2$. Let $x$ be an element of the head. The inclusion maximum element of the family that contains $A_1$ and is disjoint from $x$ must be $S_1:=A_1\cup\{y,z\}$. Repeating the construction replacing $x$ by $y$ we get $S_1'=A_1\cup\{p,q\}$. The union of $S_1$ and $S_1'$ must be in the $G'_\delta$-family. As it has cardinality strictly greater than $\delta+1$, it must be the complement of $A_2$. Hence, the head has cardinality at most four: $\delta+2\leq 4$, which is a contradiction.

Summing up, for $\delta>2$ there are only two $G'_\delta$-families with corank $3$. One is for $\delta=4$ on the ground set of cardinality $12$ subdivided into four $3$-ears. The corresponding matroid is the graphic matroid of four paths of length three that join two vertices. Good $3$-ear contractions give matroids on ground sets of cardinalities $10$, $8$, $6$, and $4$.
	
The second one is for $\delta=3$ on the set of cardinality $9$ where one element is in the head, and the others form two blocks of four elements each, each block subdivided into two $2$-ears. The corresponding matroid is graphic: two vertices $a,b$ are joined by an edge that is the head. Additionally, there are two paths of length two ($2$-ears) from $a$ to the vertex $c$ and two paths of length two from $b$ to the vertex $c$. Good $2$-ear contractions give matroids on ground sets of cardinalities $8$, $7$ (two different matroids depending if two $2$-ears were chosen from the same block or no), $6$, and $5$.
	
By the discussion from the beginning of the section the classification is complete.
\end{proof}


\end{document}